\documentclass[a4paper,12pt,reqno]{amsart}

\usepackage{amsfonts}
\usepackage{amsmath}
\usepackage{amssymb}
\usepackage{cite}
\usepackage{mathrsfs}
\usepackage{enumitem}
\usepackage[colorlinks]{hyperref}

\numberwithin{equation}{section}

\usepackage{graphicx}

\newtheorem{theorem}{Theorem}[section]

\newtheorem{remark}[theorem]{Remark}

\newtheorem{prop}[theorem]{Proposition}
\newtheorem{lemma}[theorem]{Lemma}
\newtheorem{example}[theorem]{Example}

	\allowdisplaybreaks
\begin{document}

\title[Space-time fractional diffusion equation]
{Global well-posedness of  Space-time fractional diffusion equation with Rockland operator on Graded Lie Group}
\author[A. Dasgupta]{Aparajita Dasgupta}
\address{
	Aparajita Dasgupta:
	\endgraf
	Department of Mathematics
	\endgraf
	Indian Institute of Technology, Delhi, Hauz Khas
	\endgraf
	New Delhi-110016 
	\endgraf
	India
	\endgraf
	{\it E-mail address} {\rm adasgupta@maths.iitd.ac.in}
}
\author[M. Ruzhansky]{Michael Ruzhansky}
\address{
	Michael Ruzhansky:
	\endgraf
	Department of Mathematics: Analysis, Logic and Discrete Mathematics
	\endgraf
	Ghent University, Belgium
	\endgraf
	and
	\endgraf
	School of Mathematical Sciences
	\endgraf
	Queen Mary University of London
	\endgraf
	United Kingdom
	\endgraf
	{\it E-mail address} {\rm ruzhansky@gmail.com}
}
\author[A. Tushir]{Abhilash Tushir}
\address{
	Abhilash Tushir:
	\endgraf
	Department of Mathematics
	\endgraf
	Indian Institute of Technology, Delhi, Hauz Khas
	\endgraf
	New Delhi-110016 
	\endgraf
	India
	\endgraf
	{\it E-mail address} {\rm abhilash2296@gmail.com}
}
\begin{abstract}
In this article, we examine the general space-time fractional diffusion equation for left-invariant hypoelliptic homogeneous operators on graded Lie groups. Our study covers important examples such as the time-fractional diffusion equation, the space-time fractional diffusion equation when diffusion is under the influence of sub-Laplacian on the Heisenberg group, or general stratified Lie groups. We establish the global well-posedness of the Cauchy problem for the general space-time fractional diffusion equation for the Rockland operator on a graded Lie group in the associated Sobolev spaces. More precisely,  we establish the existence and uniqueness results for both homogeneous and inhomogeneous fractional diffusion equations. In addition, we also develop some regularity estimates. 
 \end{abstract}
 
\maketitle
\tableofcontents
\section{Introduction}\label{intro}
The general Caputo-type fractional derivative $\mathbb{D}_{(g)}$, first introduced by Anatoly N. Kochubei in \cite{2011Koch}, is  defined as
\begin{equation}\label{genfracdef}
	\mathbb{D}_{(g)} u(t)=\frac{d}{d t} \int_0^t g(t-\tau) u(\tau) \mathrm{d} \tau-g(t) u(0)=\int_{0}^{t}g(t-\tau)u^{\prime}(\tau)\mathrm{d}\tau,\quad u\in C^{1}([0,\infty)),
	\end{equation}  
where the kernel function $g\in L^{1}_{loc}(0,\infty)$ satisfies  the following conditions:
\begin{enumerate}
	\item\label{c1} the Laplace transform $\widetilde{g}(p)=\mathscr{L}(g ; p)$ of $g$ 
	exists for all $p>0$;
	\item\label{c2} $\widetilde{g}(p)$ is a Stieltjes function;
	\item\label{c3} $\widetilde{g}(p) \rightarrow 0$ and $p \widetilde{g}(p) \rightarrow \infty$ as $p \rightarrow \infty$; and
	\item\label{c4} $\widetilde{g}(p) \rightarrow \infty$ and $p \widetilde{g}(p) \rightarrow 0$ as $p \rightarrow 0$.
\end{enumerate}   
If we choose the kernel function $$g(t)=\frac{t^{-\alpha}}{\Gamma(1-\alpha)},\quad 0<\alpha<1,$$ the operator $\mathbb{D}_{(g)}$ coincides with the Caputo fractional differential operator ${^{\rm{C}}D}^{\alpha}_{t}$. For more details about the general Caputo-type fractional derivative, one can refer to \cite{2011Koch}. 

\sloppy The general Caputo-type fractional derivative generalizes multi-term and distributed order fractional differential operators. There are several variations of general Caputo-type fractional derivatives, reliant upon the choice of kernel function $g$, such as the Caputo-Dzhrbashyan derivative, Caputo-Fabrizio derivative,  etc.; see \cite[Section 1]{diet20}. 
These generalized differential operators have been utilized for relaxation equations, diffusion equations, evolution equations, etc. to represent  several models of physics and phenomena, such as fractals, optics, finance, fluid mechanics, and signal processing; see \cite{chung2018,chung2019,chung2020,Ourpaper:Genfrac,MAINARDI,2023niyaz2} and references therein. We refer to Section \ref{sec:prelim}, for more details concerning the general fractional differential equations

The key objective of this paper is to study the global well-posedness of the general space-time  fractional diffusion equation for the positive Rockland operators $\mathcal{R}$ (a left-invariant hypoelliptic differential operator of homogeneous degree $\nu$) on graded Lie group $\mathbb{G}$; that is, for $s,T>0$, we study the following Cauchy problem:
\begin{equation}\label{mainpde}
	\left\{\begin{array}{l}
		\mathbb{D}_{(g)} u(t, x)+a(t)\mathcal{R}^{s} u(t, x)+b(t)u(t,x)
		=0, \quad(t, x) \in(0, T] \times  \mathbb{G}, \\
		u(0, x)=u_{0}(x), \quad x \in  \mathbb{G},
	\end{array}\right.
\end{equation}
where  $a = a(t) \geq 0$ is the diffusion coefficient and $b(t)\geq 0$ is an external potential.  Positive time-dependent potential $b(t)$ as a multiplicative factor enhances the effective diffusion rate. Modeling systems with temporally fluctuating external fields, such as time-varying magnetic or electric fields in plasmas, becomes simpler by the Cauchy problem \eqref{mainpde}.  This study extends the diffusion equations beyond conventional Euclidean space by putting the analysis within the framework of graded Lie groups.  Indeed, there are a number of generalizations and extensions of classical PDEs such as the wave equation, Klein-Gordon equation, heat equation, and Schrödinger equations in the context of graded Lie groups; see \cite{graded:chiara,graded:marianna,graded:marianna2,graded;nurgissa}. In fact, by introducing the Rockland operator and the framework of graded Lie groups, this work marks an advancement and generalization in the study of diffusion equations. It should be highlighted that the Cauchy problem \eqref{mainpde}:
\begin{itemize}
		\item for $\mathbb{G}=\mathbb{R}$,  $a\equiv 0$,  and $b(t)\equiv \lambda>0$ was the setting of \cite{2011Koch}; and\medskip
	\item for $\mathbb{G}=\mathbb{R}^{n}$, $\mathcal{R}=-\Delta$, $s=1$, $a(t)\equiv 1$,  and $b(t)\equiv 0$ was the setting of \cite{chung2020}.\medskip
\end{itemize}
Additionally, the estimates of \cite{chung2020} are included in this article as a special case. The Cauchy problem \eqref{mainpde} may also provide more general diffusion equations with adjustments to the graded Lie group and Rockland operator, for example, fractional Laplacian, sub-Laplacian on stratified Lie group, etc. In Section \ref{sec:prelim}, we provide a thorough synopsis for variants of the well-known graded Lie group and Rockland operator.
 Let us now explore the possibility of involving external forces as a source term in a fractional diffusion equation. The diffusion process may be considerably impacted when a source term is included in the fractional diffusion equation. These impacts include new equilibrium states, inhomogeneous diffusion patterns, temporal fluctuations, scaling property changes, and changes in localization. Although adding a source term brings a significant amount of complexity to the fractional diffusion equation, it also makes it possible to represent diffusion processes more accurately in the real world. The non-homogeneous general fractional diffusion equation with the Rockland operator on a graded Lie group is given by the Cauchy problem:
\begin{equation}\label{mainpdenonhom}
	\left\{\begin{array}{l}
		\mathbb{D}_{(g)} u(t, x)+a(t)\mathcal{R}^{s} u(t, x)+b(t)u(t,x)
		=f(t,x), \quad(t, x) \in(0, T] \times  \mathbb{G}, \\
		u(0, x)=u_{0}(x), \quad x \in  \mathbb{G},
	\end{array}\right.
\end{equation}
where  $a = a(t) \geq 0$ is the diffusion coefficient, $b(t)\geq 0$ is an external potential, and $f$ is the source term. It is important to note that the Cauchy problem \eqref{mainpde} for $\mathbb{D}_{(g)}$ being the  multi-term fractional differential operator, $a\equiv 1$ and $b\equiv 0$ was considered in \cite{Ruzn2020}. 

However, the hypoelliptic operators may also be taken into consideration in the general settings  of the Cauchy problems \eqref{mainpde} and \eqref{mainpdenonhom}. More precisely, we obtain the results in the following contexts:
\begin{enumerate}
	\item\label{a1} $\mathcal{R}$ being a positive Kohn-Laplacian/sub-Laplacian on the Heisenberg group $\mathbb{G}=\mathbb{H}^{n_{o}}$;
	\item\label{a2} $\mathcal{R}$ is a positive sub-Laplacian on a stratified Lie group $\mathbb{G}$; and
	\item\label{a3} $\mathcal{R}$ is any positive Rockland operator on a graded Lie group $\mathbb{G}$.
\end{enumerate}
Actually, our investigations are for the last case \eqref{a3}; the cases \eqref{a1} and \eqref{a2} appear as special cases.
\begin{remark}
 Here, we are considering  homogeneous and non-homogeneous fractional diffusion equations separately so that we can cover a wide range of diffusion equations. 
Moreover, in Section \ref{sec:results}, we will see that we can work on the homogeneous diffusion equation \eqref{mainpde} with slightly relaxed assumptions on coefficients than the non-homogeneous diffusion equation \eqref{mainpdenonhom}. Later, in Section \ref{remark}, we will highlight under what circumstances one can relax these assumptions for the non-homogeneous case.
\end{remark}
To present our main results, some background information concerning the graded Lie group and Rockland operator is required, which will be addressed in the following section. Then, we will highlight our key findings.

 We wrap up the introduction with a brief overview of the paper's structure.
\begin{itemize}
\item In Section \ref{sec:graded}, we recall the notions of graded Lie group, Rockland operator, and the associated basic Fourier analysis. Additionally, we also recall an important class of function spaces named $\mathcal{R}$-Sobolev spaces.
\item 	In Section \ref{sec:results}, we provide our main results concerning the homogeneous and non-homogeneous general fractional diffusion equations \eqref{mainpde} and \eqref{mainpdenonhom}, respectively.
	\item  In Section \ref{sec:prelim}, we review several essential tools and significant findings concerning the general Caputo-type fractional derivative. 
	\item In Section \ref{sec:main results}, we supply the proof of our main results.  Specifically, we prove the uniqueness and existence of the Cauchy problems \eqref{mainpde} and \eqref{mainpdenonhom}, and additionally, we derive certain regularity estimates for the same.
	\item Finally, in Section \ref{remark}, we make a few comments regarding the difference in assumptions on coefficients in the homogeneous and non-homogeneous homogeneous Cauchy problems.
\end{itemize}
\section{Fourier analysis of graded Lie groups and Rockland operators}\label{sec:graded}
In this section, we will review several key tools  related to graded Lie groups and  Rockland operators. A comprehensive review of the graded Lie group and Rockland operator has been explored by Folland in \cite{Folland82} and by Fischer and Ruzhansky in \cite{quantnil}. We are extracting the essential information from it, to make the paper self-contained.
 
  Let $\mathfrak{g}$ be a Lie algebra of a  connected and simply connected Lie group $\mathbb{G}$. We say $\mathbb{G}$ is a graded Lie group if its Lie algebra $\mathfrak{g}$ has a vector space decomposition
\begin{equation}\label{grad}
	\mathfrak{g}=\bigoplus\limits_{j=1}^{\infty} \mathfrak{g}_{j},\text{ such that } [\mathfrak{g}_{i},\mathfrak{g}_{j}]\subset \mathfrak{g}_{i+j},
\end{equation}
with all but finitely many of $\mathfrak{g}_{j}$ are equal to $\{0\}.$ In \cite{Folland74} and \cite{Folland75}, Folland and Stein analyzed a special case of graded Lie group called  stratified Lie group: if $\mathfrak{g}_{1}$ itself generates $\mathfrak{g}$ as an algebra, then the group $\mathbb{G}$ is said to be a stratified Lie group. Before discussing the other structures on the graded Lie group, let us look at some examples of graded Lie groups:
\begin{example}\label{expgrad}
	The Abelian group $(\mathbb{R}^{n},+)$, the Heisenberg group $\mathbb{H}^{n_{o}}$,   the group of $n \times n$ matrices that are upper triangular with $1$ on the diagonal,  the Engel group $\mathscr{B}_{4}$, etc. are some well-known examples of  graded Lie groups.
\end{example}
  The structure \eqref{grad} yields a family of dilations $\{D_{r}:\mathfrak{g}\to\mathfrak{g},r>0\}$ on graded Lie algebra $\mathfrak{g}$; see \cite[Definition 3.1.7]{quantnil}. Additionally,  graded Lie algebras  incorporate dilations naturally: if the Lie algebra $\mathfrak{g}$ is equipped with gradation  
\eqref{grad}, then we define the dilations by
\begin{equation*}
	D_{r}:=\operatorname{Exp}(A\operatorname{ln}r)=\sum\limits_{n=0}^{\infty}\frac{1}{n!}\left(\operatorname{ln}(r)A\right)^{n},\quad r>0,
\end{equation*}
where the operator $A$ is given by $AX=jX$ for $X\in\mathfrak{g}_{j}$. If the Lie algebra of a connected simply connected Lie group is equipped with dilations, then the Lie group is said to be homogeneous. Furthermore, it's also true that every graded Lie group is homogeneous as well as nilpotent, and consequently, in the light of \cite[Proposition 1.6.6]{quantnil}, the exponential map
\begin{equation*}
	\operatorname{exp}_{\mathbb{G}}:\mathfrak{g}\to\mathbb{G}\text{ defined by } x:= \operatorname{exp}_{\mathbb{G}}(x_{1}X_{1}+\cdots+x_{n}X_{n}),
\end{equation*}
is a diffeomorphism, where $\{X_{1},\dots,X_{n}\}$ is a basis of $\mathfrak{g}$. Using the above exponential map, we may transfer the dilations from the Lie algebra $\mathfrak{g}$ to the group $ \mathbb{G}$ in the following manner: the maps
\begin{equation*}
	\operatorname{exp}_{\mathbb{G}}\circ D_{r}\circ\operatorname{exp}_{\mathbb{G}}^{-1}:\mathbb{G}\to\mathfrak{g}\to\mathfrak{g}\to \mathbb{G},\quad r>0,
\end{equation*}
are automorphisms of the group $\mathbb{G}$. We will also refer to them as dilations on $\mathbb{G}$ and denote them by $D_r$. We can write
\begin{equation*}
	D_{r}(x):=rx=(r^{\nu_{1}}x_{1},\dots,r^{\nu_{n}}x_{n}),\quad x\in\mathbb{G},
\end{equation*}
where $\nu_{1},\dots,\nu_{n}$ are the weights of the dilations. 

To clarify the aforementioned abstract concepts, let us use a few specific examples from Example \ref{expgrad}:
\begin{example} The abelian group $(\mathbb{R}^{n},+)$ is a graded Lie group as its Lie algebra $\mathfrak{g}=\mathbb{R}^{n}$ is trivially graded, i.e., $\mathfrak{g}_{1}=\mathbb{R}^{n}$ and $\mathfrak{g}_{j}=\{0\}$ for all $j\geq 2$. Thus, $(\mathbb{R}^{n},+)$ is a stratified Lie group. Furthermore, it is a homogeneous Lie group when equipped with the usual dilations
	\begin{equation*}
	 D_{r}x = rx, r > 0,\quad x\in\mathbb{R}^{n}.
	\end{equation*}
\end{example}
\begin{example}
	Let $n_{o}\in\mathbb{N}$ and $\mathbb{H}^{n_{o}}$ be the Heisenberg group whose underlined space is $\mathbb{R}^{2n_{o}+1}\sim\mathbb{R}^{n_{o}}\times\mathbb{R}^{n_{o}}\times\mathbb{R}$ equipped with the group law
	 \begin{equation*}
		(x, y, t)\left(x^{\prime}, y^{\prime}, t^{\prime}\right):=\left(x+x^{\prime}, y+y^{\prime}, t+t^{\prime}+\frac{1}{2}\left(x y^{\prime}-x^{\prime} y\right)\right),
	\end{equation*}
	where $(x, y, t)$, $\left(x^{\prime}, y^{\prime}, t^{\prime}\right)\in\mathbb{R}^{n_{o}} \times \mathbb{R}^{n_{o}} \times \mathbb{R} \sim \mathbb{H}^{n_{o}}$. The Lie algebra $\mathfrak{h}_{n_{o}}$ of the Heisenberg group is $\mathbb{R}^{2n_{o}+1}$ with a canonical basis
	\begin{equation}\label{basis}
		\mathfrak{h}_{n_{o}}=\{X_{1},\dots,X_{n},Y_{1},\dots,Y_{n},T\},
	\end{equation} 
	satisfying the commutator relation
	\begin{equation*}
		[X_{j},Y_{j}]=T,\quad j=1,\dots,n_{o},
	\end{equation*}
\sloppy	and all the other Lie brackets  are trivial (excluding those derived through anti-symmetry). Further, the Heisenberg group $\mathbb{H}^{n_{o}}$ is graded as its Lie algebra $\mathfrak{h}_{n_{o}}$ admits the following vector space decomposition:
\begin{equation*}
	\mathfrak{h}_{n_{o}}=\mathfrak{g}_{1}\oplus\mathfrak{g}_{2}=\operatorname{span}\{X_{1},\dots,X_{n},Y_{1},\dots,Y_{n}\}\oplus
	\mathbb{R}T.
\end{equation*}
The dilations of Lie algebra  $\mathfrak{h}_{n_{o}}$, i.e., $D_{r}:\mathfrak{h}_{n_{o}}\to \mathfrak{h}_{n_{o}}$ are given by
\begin{equation*}
	D_{r}(X_{j}) = rX_{j},\quad D_{r}(Y_{j}) = rY_{j},\quad  j = 1,\dots,n_{o}, \text{ and } D_{r}(T) = r^{2}T.    
\end{equation*}
Consequently, the dilations of the Heisenberg group $\mathbb{H}^{n_{o}}$, i.e., $D_{r}:\mathbb{H}^{n_{o}}\to \mathbb{H}^{n_{o}}$, are given by
\begin{equation*}
	D_{r}(h)=rh=(rx,ry,r^{2}t),\quad (x,y,t)\in\mathbb{R}^{n_{o}}\times\mathbb{R}^{n_{o}}\times \mathbb{R}.
\end{equation*}
\end{example}

\medskip

Let us now discuss the Rockland operator, a natural differential operator on a graded Lie group, in detail.
 Let $\widehat{\mathbb{G}}$ be the unitary dual of $\mathbb{G}$ and $\mathcal{H}^{\infty}_{\pi}$ be the space of smooth vectors; see \cite[Definition 1.7.2]{quantnil}. A Rockland operator $\mathcal{R}$ on $\mathbb{G}$ is a left-invariant, homogeneous differential operator of positive degree $ \nu$ that meets the following Rockland condition:
\begin{equation*}
\text{for each } \pi\in \widehat{\mathbb{G}}\setminus\{1\} \text{ and }	\forall v\in\mathcal{H}^{\infty}_{\pi},\quad \pi(\mathcal{R})v=0\implies v=0.
\end{equation*}
There have been alternative characterizations of these operators studied by Beals in \cite{Beals} and Rockland in \cite{rockland87}. For general nilpotent Lie groups, the existence of Rockland operators perfectly characterizes the class of graded Lie groups; see \cite{Miller80} and \cite{Robinson97}. Let us first examine a few examples of Rockland operators before continuing on to the other properties of Rockland operators. 
	\begin{itemize}
	\item For $\mathbb{G}=(\mathbb{R}^{n},+)$, let $P$ be a homogeneous polynomial that does not vanish except at zero. Then
	\begin{equation*}
		P(-i\partial_{1},\dots,-i\partial_{n}) \text{ is a Rockland operator.}
	\end{equation*}
	For example, we can take $\mathcal{R}=(-\Delta)^{m}$ or $\mathcal{R}=(-1)^{m}\sum\limits_{j=1}^{n}\left(\frac{\partial}{\partial x_{j}}\right)^{2m}$.
	\item Any sub-Laplacian on a stratified Lie group is a Rockland operator of homogeneous degree $2$.
	\item For $\mathbb{G}=\mathbb{H}^{n_{o}}$, we can take
	\begin{equation*}
		\mathcal{R}=(-\mathcal{L})^{m} \text{ or } \mathcal{R}=(-1)^{m}\sum\limits_{j=1}^{n}\left(X^{2m}_{j}+Y_{j}^{2m}\right).
	\end{equation*}
	where $\mathcal{L}$ is the sub-Laplacian and $X_{j},Y_{j}$ are given in \eqref{basis}.
	\item For $k\in\mathbb{N}$, the operators $\mathcal{R}^{k}$ and $\overline{\mathcal{R}}$ are also Rockland operators.
	\end{itemize}
According to \cite[Proposition 4.1.15]{quantnil}, recall that the Rockland operator $\mathcal{R}$ is densely defined on $\mathcal{D}(\mathbb{G})\subset L^{2}(\mathbb{G})$ while the operator $\pi(\mathcal{R})$ is densely defined on $\mathcal{H}^{\infty}_{\pi}\subset \mathcal{H}_{\pi}$. Taking into account \cite[Theorem VIII.6]{Simon80}, i.e., the spectral theorem for unbouded operators, we can write
\begin{equation*}
	\mathcal{R}=\int_{\mathbb{R}}\lambda\mathrm{d}E(\lambda)\text{ and } 	\pi(\mathcal{R})=\int_{\mathbb{R}}\lambda\mathrm{d}E_{\pi}(\lambda),
\end{equation*}
where $E$ and $E_{\pi}$ are the spectral measures linked to $\mathcal{R}$ and $\pi(\mathcal{R})$.  Hulanicki et. al. in \cite{Hulan85} investigated the spectrum of $\pi(\mathcal{R})$ and proved that it is purely discrete and positive for $\pi\in\widehat{\mathbb{G}}\setminus\{1\}.$ Moreover, it admits the following infinite matrix representation:
 \begin{equation}\label{piofrs1}
\pi(\mathcal{R})=\operatorname{diag}\left(\pi_{1}^{2},\pi_{2}^{2},\dots\right),
\end{equation}
consequently, we have
 \begin{equation}\label{piofrs}
	\pi(\mathcal{R}^{s})=\operatorname{diag}\left(\pi_{1}^{2s},\pi_{2}^{2s},\dots\right),\quad s>0.
\end{equation}

\medskip

Let us now discuss the basic Fourier analysis related to the graded Lie group and Rockland operator.
Let $f\in L^{1}(\mathbb{G})$ and $\pi\in\widehat{\mathbb{G}}$, the group Fourier transform of $f$ at $\pi$
is defined by
\begin{equation}\label{grpft}
	\mathcal{F}_{\mathbb{G}}f(\pi)\equiv \widehat{f}(\pi)\equiv  \pi(f):=\int\limits_{\mathbb{G}}f(x)\pi(x)^{*}\mathrm{d}x,
\end{equation}
 where integration is performed against the bi-invariant Haar measure on $\mathbb{G}$. The linear mappling $\pi(f):\mathcal{H}_{\pi}\to \mathcal{H}_{\pi}$ may be written as an infinite matrix as soon as we specify a basis for Hilbert space $\mathcal{H}_{\pi}$. Furthermore, for any $f\in \mathcal{S}(\mathbb{G})\cap L^{1}(\mathbb{G})$ and $X\in\mathfrak{g}$, we obtain
 \begin{equation*}
 	\mathcal{F}_{\mathbb{G}}\left(Xf\right)(\pi)=\pi(X)\pi(f)=\pi(X)\widehat{f}(\pi).
 \end{equation*}
 As a consequence, we have
 \begin{equation*}
 	\mathcal{F}_{\mathbb{G}}\left(\mathcal{R}f\right)(\pi)=\pi(\mathcal{R})\widehat{f}(\pi)\text{ and } 	\mathcal{F}_{\mathbb{G}}\left(\mathcal{R}^{s}f\right)(\pi)=\pi(\mathcal{R}^{s})\widehat{f}(\pi),
 \end{equation*}
 and  combining it with the relations \eqref{piofrs1} and \eqref{piofrs}, we obtain
 \begin{equation}\label{ftrock}
	\mathcal{F}_{\mathbb{G}}\left(\mathcal{R}f\right)(\pi)=\left\{\pi_{m}^{2}\cdot\widehat{f}(\pi)_{m,k}\right\}_{m,k\in\mathbb{N}}\text{ and } 	\mathcal{F}_{\mathbb{G}}\left(\mathcal{R}^{s}f\right)(\pi)=\left\{\pi_{m}^{2s}\cdot\widehat{f}(\pi)_{m,k}\right\}_{m,k\in\mathbb{N}}.
 \end{equation}
  The Fourier inversion formula is given by
 \begin{equation}\label{fourierinversion}
 	f(x)=\int_{\widehat{\mathbb{G}}}\operatorname{Tr}\left[\pi(x)\widehat{f}(\pi)\right]\mathrm{d}\mu(\pi),
 \end{equation}
where $\mathrm{d}\mu$ is the Plancherel measure on $\widehat{\mathbb{G}}.$ Additionally,  the linear mapping $\pi(f):\mathcal{H}_{\pi}\to \mathcal{H}_{\pi}$ is a Hilbert-Schmidt operator:
\begin{equation*}
	\|\pi(f)\|^{2}_{\mathrm{HS}(\mathcal{H}_{\pi})}=\operatorname{Tr}\left[\pi(f)\pi(f)^{*}\right]<\infty.
\end{equation*}
 And also, the following Plancherel formula holds:
  \begin{equation}\label{planch}
  	\int_{\mathbb{G}}|f(x)|^{2}\mathrm{d}x=\int_{\widehat{\mathbb{G}}}\|\pi(f)\|_{\operatorname{HS}}^{2}\mathrm{d}\mu(\pi).
  \end{equation}
Next, we must review the notion of Sobolev spaces associated with the Rockland operator, adapted for the graded Lie group. Thus, let $\mathcal{R}$ be a positive homogeneous Rockland operator of degree $\nu$. For $\gamma\in\mathbb{R}$, the $\mathcal{R}$-Sobolev space $L^{2}_{\gamma}(\mathbb{G})$ is identified as the space of tempered distribution $\mathcal{S}^{\prime}(\mathbb{G})$ obtained by the completion of the Schwartz space $\mathcal{S}(\mathbb{G})$ w.r.t. the following Sobolev norm:
\begin{equation}\label{sobnorm}
	\|f\|_{L^{2}_{\gamma}(\mathbb{G})}:=\left\|(I+\mathcal{R})^{\frac{\gamma}{\nu}}f\right\|_{L^{2}(\mathbb{G})},\quad \gamma\in\mathbb{R}.
\end{equation}
Moreover, taking into account \cite[Theorem 4.4.3]{quantnil},  for $\gamma\geq0$, the following norms are equivalent to $\|\cdot\|_{L^{2}_{\gamma}(\mathbb{G})}$:
\begin{equation}\label{eqult}
	f\mapsto \|f\|_{L^{2}(\mathbb{G})}+\left\|(I+\mathcal{R})^{\frac{\gamma}{\nu}}f\right\|_{L^{2}(\mathbb{G})},\quad f\mapsto \|f\|_{L^{2}(\mathbb{G})}+\left\|\mathcal{R}^{\frac{\gamma}{\nu}}f\right\|_{L^{2}(\mathbb{G})}.
\end{equation}
In the light of \cite[Corollary 4.3.11 (ii)]{quantnil}, for $\gamma\geq0$, we also deduce that
\begin{equation*}
	\|f\|_{L^{p}(\mathbb{G})}\leq C\|f\|_{L^{p}_{\gamma}(\mathbb{G})},\quad f\in\mathcal{S}(\mathbb{G}).
\end{equation*}
 Moreover, the homogeneous/inhomogeneous $L^{p}$-Sobolev spaces on graded Lie group $\mathbb{G}$ with any Rockland operators coincide. In other words, the Sobolev spaces do not depend on the Rockland operator $\mathcal{R}$; see \cite[Theorem 4.4.20]{quantnil}.
  \section{Main results}\label{sec:results}
  In this section, we establish the global well-posedness of the Cauchy problem for the general space-time fractional diffusion equation for the Rockland operator on a graded Lie group in the associated $\mathcal{R}$-Sobolev spaces. More precisely, we establish the existence and uniqueness of the generalized solution to the Cauchy problems \eqref{mainpde} and \eqref{mainpdenonhom}. Moreover, we also develop some regularity estimates for the solutions.  We shall write $A\lesssim B$ throughout the paper if there exists a positive constant $C$ that is independent of the  parameters that occur such that $A \leq CB$.  
  
The following theorem gives the well-posedness result for the homogeneous Cauchy problem \eqref{mainpde}:
\begin{theorem}[Homogeneous case]\label{mthm1}
	Let $\mathbb{G}$ be a graded Lie group	and $\mathcal{R}$ be a Rockland operator of homogeneous degree $\nu$. Assume that the initial Cauchy data satisfies $u_{0}\in L^{2}(\mathbb{G})$ and  also assume that $0 \leq a(t),b(t)\in C([0,T])$   such that  either $\inf\limits_{t\in[0,T]}a(t)=a_{0}>0$ or $\inf\limits_{t\in[0,T]}b(t)=b_{0}>0$.  Then the Cauchy problem \eqref{mainpde} has a unique solution $u\in C([0,T];L^{2}(\mathbb{G}))$. Furthermore, we have the following regularity estimates:
	
	\begin{enumerate}
		\item if $u_{0}\in L^{2}(\mathbb{G})$, then the solution $u$  satisfies the
		estimate  
		\begin{equation}\label{MI01}
			\|u(t,\cdot)\|_{L^{2}(\mathbb{G})}\lesssim \|u_{0}\|_{L^{2}(\mathbb{G})};
		\end{equation}
		\item if $u_{0}\in L^{2}_{\gamma}(\mathbb{G})$ for some $\gamma\in\mathbb{R}$, then the solution $u$  satisfies the
		estimate  
		\begin{equation}\label{MI02}
			\|u(t,\cdot)\|_{L^{2}_{\gamma}(\mathbb{G})}\lesssim \|u_{0}\|_{L^{2}_{\gamma}(\mathbb{G})};
		\end{equation}
		\item if $u_{0}\in L^{2}_{s\nu}(\mathbb{G})$, then the solution $u$  satisfies the
		estimate  
		\begin{equation}\label{MI03}
			\|\mathbb{D}_{(g)}u(t,\cdot)\|_{L^{2}(\mathbb{G})}\lesssim \|u_{0}\|_{L^{2}_{s\nu}(\mathbb{G})};
		\end{equation}
		\item if $u_{0}\in  L^{2}_{\gamma+s\nu}(\mathbb{G})$ for some $\gamma\in\mathbb{R}$, then the solution $u$  satisfies the
		estimate  
		\begin{equation}\label{MI04}
			\|\mathbb{D}_{(g)}u(t,\cdot)\|_{L^{2}_{\gamma}(\mathbb{G})}\lesssim \|u_{0}\|_{L^{2}_{\gamma+s\nu}(\mathbb{G})};
		\end{equation}
	\end{enumerate}
	for all $t\in[0,T]$.
\end{theorem}
The following theorem gives the well-posedness result for the non-homogeneous Cauchy problem \eqref{mainpdenonhom}:
\begin{theorem}[Non-homogeneous case]\label{mthm2}
	Let $\mathbb{G}$ be a graded Lie group	and $\mathcal{R}$ be a Rockland operator of homogeneous degree $\nu$. Assume that the initial Cauchy data satisfies $u_{0}\in L^{2}(\mathbb{G})$ and the source term $f\in C([0,T];L^{2}(\mathbb{G}))$. If we  also assume that $0\leq a(t),b(t)\in C([0,T])$ such that  $\inf\limits_{t\in[0,T]}b(t)=b_{0}>0$,  then the Cauchy problem \eqref{mainpdenonhom} has a unique  solution $u\in C([0,T];L^{2}(\mathbb{G}))$. Furthermore, we have the following regularity estimates:
	
		\begin{enumerate}
		\item if $u_{0}\in L^{2}(\mathbb{G})$ and $f\in C([0,T];L^{2}(\mathbb{G}))$, then the solution $u$  satisfies the
		estimate  
		\begin{equation}\label{I01}
\|u(t,\cdot)\|_{L^{2}(\mathbb{G})}\lesssim \|u_{0}\|_{L^{2}(\mathbb{G})}+\|f\|_{C([0,T];L^{2}(\mathbb{G}))};
		\end{equation}
		\item if $u_{0}\in L^{2}_{\gamma}(\mathbb{G})$ and $f\in C([0,T];L^{2}_{\gamma}(\mathbb{G}))$, for some $\gamma\in\mathbb{R}$, then the solution $u$  satisfies the
		estimate  
		\begin{equation}\label{I02}
			\|u(t,\cdot)\|_{L^{2}_{\gamma}(\mathbb{G})}\lesssim \|u_{0}\|_{L^{2}_{\gamma}(\mathbb{G})}+\|f\|_{C([0,T];L^{2}_{\gamma}(\mathbb{G}))};
		\end{equation}
		\item if $u_{0}\in L^{2}_{s\nu}(\mathbb{G})$ and $f\in C([0,T];L^{2}_{s\nu}(\mathbb{G}))$, then the solution $u$  satisfies the
		estimate  
		\begin{equation}\label{I03}
			\|\mathbb{D}_{(g)}u(t,\cdot)\|_{L^{2}(\mathbb{G})}\lesssim \|u_{0}\|_{L^{2}_{s\nu}(\mathbb{G})}+\|f\|_{C([0,T];L^{2}_{s\nu}(\mathbb{G}))};
		\end{equation}
		\item if $u_{0}\in L^{2}_{\gamma+s\nu}(\mathbb{G})$ and $f\in C([0,T];L^{2}_{\gamma+s\nu}(\mathbb{G}))$, for some $\gamma\in\mathbb{R}$, then the solution $u$  satisfies the
		estimate  
		\begin{equation}\label{I04}
			\|\mathbb{D}_{(g)}u(t,\cdot)\|_{L^{2}_{\gamma}(\mathbb{G})}\lesssim \|u_{0}\|_{L^{2}_{\gamma+s\nu}(\mathbb{G})}+\|f\|_{C([0,T];L^{2}_{\gamma+s\nu}(\mathbb{G}))};
		\end{equation}
	\end{enumerate}
	for all $t\in[0,T]$.
\end{theorem}
\noindent
Estimates in Theorem \ref{mthm1} may be derived from Theorem \ref{mthm2} when $\inf\limits_{t\in[0,T]}b(t)=b_{0}>0$, but not when $b(t)\geq 0$. In Section \ref{remark}, we shall discuss certain circumstances when estimates in Theorem \ref{mthm1} can be recovered from Theorem \ref{mthm2} for the case $b(t)\geq 0$.
\section{Preliminaries}\label{sec:prelim} 
In this section, we review several essential tools and significant findings concerning the general Caputo-type fractional derivative.  The following is the compressed version of our recent work; for details, one can refer to \cite[Section 3]{Ourpaper:Genfrac}.

Consider the following fractional differential equation associated with the general Caputo-type fractional derivative $\mathbb{D}_{(g)}$:
\begin{equation}\label{homchy}
	\left\{\begin{array}{l}
		\mathbb{D}_{(g)} w(t)+\lambda w(t)=f(t), \quad t>0,\lambda>0, \\
		w(0)=w_{0}.
	\end{array}\right.
\end{equation}
The following lemma is dedicated to the corresponding homogeneous fractional differential equation.
\begin{lemma}\label{homthm}\cite[Theorem 2]{2011Koch}  The fractional differential equation \eqref{homchy} with initial Cauchy data $w_{0}=1$ and the source term $f\equiv 0$, 
has a  unique solution $w_{\lambda}(t)$ that satisfies the following:
\begin{enumerate}
\item\label{P1} $w_{\lambda}$ is continuous on $[0,\infty)$;
\item\label{P2} $w_{\lambda}$ is infinitely differentiable on $(0,\infty)$;  and
\item\label{P3} $(-1)^{n}w^{(n)}_{\lambda}\geq 0$  for all $t>0$ and $n\in\mathbb{N}_{0}$.
\end{enumerate}
\end{lemma}
\begin{remark}\label{firstremark}
In continuation of the above theorem, using the monotonicity of $w_{\lambda}$ and the Karamata-Feller Tauberian theorem in \cite{fel:book}, the author also proved that 
\begin{equation*}
    w_{\lambda}(t)\to 1,\quad \text{as }t\to 0^{+}.
\end{equation*}
\end{remark}
\begin{remark}\label{2ndremark}
    In light of Lemma \ref{homthm} and the above remark, it is easy to note that the solution $w_{\lambda}(t)$ satisfies
    \begin{equation*}
        0\leq w_{\lambda}(t)\leq 1, \quad \text{ for all }t>0.
    \end{equation*}
\end{remark}
\begin{lemma}\cite[Lemma 3.1]{chung2020}\label{nonthm}
	Let $\lambda>0, T>0$, and $f \in C([0, T])$. Then the fractional differential equation \eqref{homchy} with the initial Cauchy data $w_{0}=0$,
 has a unique solution in $C([0, T])$. In particular, the solution has the form
	\begin{equation*}
	w(t)=-\frac{1}{\lambda} \int_0^t  f(s)w_{\lambda}^{\prime}(t-s) \mathrm{d} s,
	\end{equation*}
	where $w_{\lambda}(t)$ is the solution to the homogeneous Cauchy problem obtained in Lemma \ref{homthm}.
\end{lemma}
\begin{remark}\label{3rdremark}
    Combining the Lemmas \ref{homthm} and \ref{nonthm}, the fractional differential equation \eqref{homchy} has a unique solution in $C([0,T])$ expressed as
    \begin{equation*}
        w(t)=w_{0}w_{\lambda}(t)-\frac{1}{\lambda} \int_0^t  f(s)w_{\lambda}^{\prime}(t-s) \mathrm{d} s.
    \end{equation*}
\end{remark}
\begin{prop}\label{proppp} Let the kernel function $g\in L^{1}_{loc}(0,\infty)$ satisfy the conditions \eqref{c1}-\eqref{c4} in the introduction and let $v\in C([0,T])$ with $\mathbb{D}_{(g)}v\in C([0,T])$. Then the following statements hold:
	\begin{enumerate}
\item If $t_{0} \in (0,T]$ and $v(t_{0}) = \max\{v(t):t\in[0,T]\}$, then we have $\mathbb{D}_{(g)}v(t_{0})\geq 0$. \medskip
\item If $t_{0} \in (0,T]$ and $v(t_{0}) = \min\{v(t):t\in[0,T]\}$, then we have $\mathbb{D}_{(g)}v(t_{0})\leq 0$.
	\end{enumerate}

\end{prop}

\section{Proof of main results}\label{sec:main results}
In this section, we first prove that the classical solution $u$ to the Cauchy problems \eqref{mainpde} and  \eqref{mainpdenonhom} exists and is unique. We next provide some required regularity estimates to wrap up the proof of Theorem \ref{mthm1} and Theorem \ref{mthm2}. 

Before proving the existence and uniqueness part, let us observe that the existence and uniqueness of Theorem \ref{mthm1} can be proved using the same justifications as Theorem \ref{mthm2} by simply placing $f\equiv0$.
Thus, we only prove the uniqueness and existence of Theorem \ref{mthm2}, and the proof of Theorem \ref{mthm1} should be considered verbatim.
\begin{proof}[Proof of Theorem \ref{mthm2}](Existence and uniqueness)
Applying the group Fourier transform \eqref{grpft} to the Cauchy problem \eqref{mainpdenonhom} with regard to $x\in\mathbb{G}$  and for all $\pi\in \widehat{\mathbb{G}}$, we can derive the set of Cauchy problems with Fourier coefficients:     
\begin{equation*}
	\left\{\begin{array}{l}
		\mathbb{D}_{(g)} \widehat{u}(t, \pi)+\pi(\mathcal{R}^{s})a(t) \widehat{u}(t, \pi)+b(t)\widehat{u}(t, \pi)=\widehat{f}(t, \pi), \quad (t,\pi)\in(0,T]\times  \widehat{\mathbb{G}}, \\
		\widehat{u}(0, \pi)=\widehat{u}_{0}(\pi), \quad \pi \in  \widehat{\mathbb{G}}.
	\end{array}\right.
\end{equation*}
With the aid of relation \eqref{ftrock} of the infintesimal representation of $\pi(\mathcal{R}^{s})$, the aforementioned Cauchy problem can be viewed concurrently as an infinite system of equations of the type
\begin{equation}\label{mkftpde}
	\left\{\begin{array}{l}
		\mathbb{D}_{(g)} \widehat{u}(t, \pi)_{m,k}+\pi^{2s}_{m}a(t) \widehat{u}(t, \pi)_{m,k}+b(t)\widehat{u}(t, \pi)_{m,k}=\widehat{f}(t, \pi)_{m,k}, \quad t\in(0,T],  \\
		\widehat{u}(0, \pi)_{m,k}=\widehat{u}_{0}(\pi)_{m,k},
	\end{array}\right.
\end{equation}
for all $\pi\in \widehat{\mathbb{G}}$ and $(m,k)\in\mathbb{N}^{2}$. At this moment, we fix an arbitrary representation $\pi\in\widehat{\mathbb{G}}$ and an ordered pair $(m,k)\in\mathbb{N}^{2}$. We then handle each of the equations provided by \eqref{mkftpde} separately.
Equivalently, we can rewrite it as
\begin{equation}\label{modtrans}
	\mathbb{D}_{(g)} \widehat{u}(t, \pi)_{m,k}+\Lambda_{\pi}^{s,m}\widehat{u}(t, \pi)_{m,k}=\widehat{h}(t, \pi)_{m,k}, \quad	\widehat{u}(0, \pi)_{m,k}=\widehat{u}_{0}(\pi)_{m,k}, \quad t\in(0,T], 
\end{equation}
where 
\begin{equation*}
	\Lambda_{\pi}^{s,m}:=\pi^{2s}_{m}a_{1}+b_{1},\quad a_{1}:=\sup\limits_{t\in[0,T]}\{a(t)\},\quad b_{1}:=\sup\limits_{t\in[0,T]}\{b(t)\},
\end{equation*}
and
\begin{equation*}
	\widehat{h}(t, \pi)_{m,k}:=\widehat{f}(t, \pi)_{m,k}+\left[\pi^{2s}_{m}\left(a_{1}-a(t)\right)+\left(b_{1}-b(t)\right)\right]\widehat{u}(t, \pi)_{m,k},
\end{equation*}
for all $(m,k)\in\mathbb{N}^{2}$. Since $\Lambda_{\pi}^{s,m}>0$, therefore, utilizing Remark \ref{3rdremark} for Cauchy problem \eqref{modtrans}, we obtain a unique generalized solution $\widehat{u}(\cdot, \pi)_{m,k}\in C([0,T])$, since $\widehat{h}(\cdot, \pi)_{m,k}\in C([0,T])$, for all $(m,k)\in\mathbb{N}^{2}$. To be more accurate, for each $\pi\in\widehat{\mathbb{G}}$ and $(m,k)\in\mathbb{N}^{2}$ , the generalized solution $\widehat{u}(t, \pi)_{m,k}$ can be explicitly written as
\begin{equation}\label{soluuxi}
	\widehat{u}(t, \pi)_{m,k}=\widehat{u}_{0}(\pi)_{m,k}\widehat{u}_{\Lambda_{\pi}^{s,m}}(t,\pi)_{m,k}-
	\frac{1}{\Lambda_{\pi}^{s,m}}\int^{t}_{0}\widehat{h}(t, \pi)_{m,k}\frac{\partial}{\partial s}\widehat{u}_{\Lambda_{\pi}^{s,m}}(t-s,\pi)_{m,k}\mathrm{d}s,
\end{equation}
where the function $\widehat{u}_{\Lambda_{\pi}^{s,m}}(t,\pi)_{m,k}$ is the solution of the  differential equation:
\begin{equation*}
		\mathbb{D}_{(g)} \widehat{u}(t, \pi)_{m,k}+\Lambda_{\pi}^{s,m}\widehat{u}(t, \pi)_{m,k}=0, \quad \widehat{u}(0, \pi)_{m,k}=1, \quad t\in(0,T], \\
\end{equation*}
obtained using Theorem \ref{homthm}.  Furthermore, we can express the solution $u(t,x)$ in the following form using the Fourier inversion formula \eqref{fourierinversion}:
\begin{equation*}
	u(t,x)=\int\limits_{\pi\in\widehat{\mathbb{G}}}\operatorname{Tr}\left[\pi(x)\widehat{u}(t,\pi)\right]\mathrm{d}\mu(\pi), \quad (t,x)\in [0,T]\times \mathbb{G},
\end{equation*}
where $\widehat{u}(t,\pi)=\{\widehat{u}(t, \pi)_{m,k}\}_{m,k\in\mathbb{N}}$ is given in \eqref{soluuxi} and $d\mu(\pi)$ is the Plancherel measure on $\widehat{\mathbb{G}}$. The uniqueness of solution $u(t,x)$ is inferred from the uniqueness of solution $\widehat{u}(t,\pi)_{m,k}$ as shown in \eqref{soluuxi}. This completes the proof.	
\end{proof}
Thus, we have proved the exitence and uniqueness of solution to the Cauchy problem \eqref{mainpdenonhom} as well as \eqref{mainpde}.

Prior to commencing the proof of regularity estimates in Theorem \ref{mthm1}, it is necessary to establish the regularity of the solution $\widehat{u}(t,\pi)_{m,k}$ using the following introductory lemma:
\begin{lemma}\label{regularity}
	Let $\pi\in\widehat{\mathbb{G}}$, $(m,k)\in\mathbb{N}^{2}$ and  $\widehat{u}(t,\pi)_{m,k}$ be the unique solution of the Cauchy problem \eqref{mkftpde}. Then the following statements hold:
	\begin{enumerate}
		\item If $\widehat{f}(t,\pi)_{m,k} \leq 0$ on $[0, T]$ and $\widehat{u}_{0}(\pi)_{m,k} \leq 0$, then we have $\widehat{u}(t,\pi)_{m,k}  \leq 0$ on $[0, T]$,
		\item If $\widehat{f}(t,\pi)_{m,k}  \geq 0$ on $[0, T]$ and $\widehat{u}_{0}(\pi)_{m,k}  \geq 0$, then we have $\widehat{u}(t,\pi)_{m,k}  \geq 0$ on $[0, T]$.
	\end{enumerate}
\end{lemma}
\begin{proof} 
	In the light of Theorem \ref{mthm1}, we have $\widehat{u}(\cdot,\pi)_{m,k}\in C([0,T])$ for all $(m,k)\in\mathbb{N}^{2}$, therefore, there exists a $t_{0}\in[0,T]$	such that
	\begin{equation*}
		\widehat{u}(t_{0},\pi)_{m,k}=\max \{\widehat{u}(t,\pi)_{m,k}:t\in[0,T]\}.
	\end{equation*}
	If $t_{0}=0$, then \begin{equation*}
		\widehat{u}(t,\pi)_{m,k}\leq \widehat{u}(0,\pi)_{m,k}=\widehat{u}_{0}(\pi)_{m,k}\leq 0.
	\end{equation*}
	From \eqref{mkftpde}, we have
	\begin{equation*}
			\mathbb{D}_{(g)} \widehat{u}(t, \pi)_{m,k}=\widehat{f}(t, \pi)_{m,k}-\pi^{2s}_{m}a(t) \widehat{u}(t, \pi)_{m,k}-b(t)\widehat{u}(t, \pi)_{m,k},\quad \text{for all }t\in[0,T].
	\end{equation*}
	Since $a,b,\widehat{f}(\cdot, \pi)_{m,k},\widehat{u}(\cdot, \pi)_{m,k}\in C([0,T])$, we have $\mathbb{D}_{(g)} \widehat{u}(t, \pi)_{m,k}\in C([0,T])$.
	If $t_{0}\in(0,T]$, then Proposition \ref{proppp} implies that $\mathbb{D}(g)\widehat{u}(t_{0},\pi)_{m,k}\geq0$, whence we get
	\begin{equation*}
		(\pi^{2s}_{m}a(t_{0})+b(t_{0})) \widehat{u}(t_{0}, \pi)_{m,k}=-\mathbb{D}_{(g)} \widehat{u}(t_{0}, \pi)_{m,k}+\widehat{f}(t_{0}, \pi)_{m,k}\leq0.
	\end{equation*}
	Since $\pi^{2s}_{m}a(t)+b(t)>0$ for all $t\in[0,T]$, we obtain $ \widehat{u}(t_{0}, \pi)_{m,k}\leq 0$. Hence 
	\begin{equation*}
		\widehat{u}(t, \pi)_{m,k}\leq \max \{\widehat{u}(t,\pi )_{m,k}:t\in[0,T]\}=\widehat{u}(t_{0},\pi )_{m,k}\leq 0,\quad \text{for all }t\in[0,T].
	\end{equation*}
	Similarily, by repeating the argument for $\widehat{u}^{*}(t,\pi)_{m,k}=-\widehat{u}(t,\pi)_{m,k}$, we can conclude the second part. This completes the proof. 
\end{proof}
We can now establish the regualarity estimates for the homogeneous Cauchy problem \eqref{mainpde}.

In order to describe the solution corresponding to the Cauchy problems with coefficients $a(t),b(t)$ and $a_{0},b_{0}$, respectively, we will specify the coefficients as $u(t,x;a,b)$ and $u(t,x;a_{0},b_{0})$ throughout the proof of  Theroem \ref{mthm1}.
\begin{proof}[Proof of Theorem \ref{mthm1}](Regularity estimates) Consider the homogeneous Cauchy problem
	\begin{equation}\label{homeqnab}
		\left\{\begin{array}{l}
			\mathbb{D}_{(g)} u(t, x;a,b)+a(t)\mathcal{R}^{s} u(t, x;a,b)+b(t)u(t, x;a,b)=0, \quad t \in(0, T],x\in \mathbb{G}, \\
			u(0, x;a,b)=u_{0}(x), \quad x \in  \mathbb{G}.
		\end{array}\right.
	\end{equation}
	Applying the group Fourier transform \eqref{grpft} and using the relation \eqref{ftrock}, we get
	\begin{equation}\label{maintranspde1}
		\left\{\begin{array}{l}
			\mathbb{D}_{(g)} \widehat{u}(t, \pi;a,b)_{m,k}+\pi_{m}^{2s}a(t) \widehat{u}(t, \pi;a,b)_{m,k}+b(t)\widehat{u}(t, \pi;a,b)_{m,k}=0,~~t\in(0,T], \\
			\widehat{u}(0, \pi;a,b)_{m,k}=\widehat{u}_{0}(\pi)_{m,k},
		\end{array}\right.
	\end{equation}
for all $(m,k)\in\mathbb{N}^{2}$. At this moment, we fix an arbitrary representation $\pi\in\widehat{\mathbb{G}}$ and an ordered pair $(m,k)\in\mathbb{N}^{2}$. We then handle each of the equations provided by \eqref{maintranspde1} separately. It is crucial to keep in mind that $\widehat{u}(t, \pi;a,b)_{m,k}$ is merely a function of $t$. Using Lemma \ref{regularity} for the Cauchy problem \eqref{maintranspde1} allows us to conclude that
	\begin{equation}\label{signreg}
		\widehat{u}(t, \pi;a,b)_{m,k} ~~\text{ holds the same sign as }~~  \widehat{u}_{0}(\pi)_{m,k}.
	\end{equation}
	For $a_{0}:=\inf\limits_{t\in[0,T]}a(t)$ and $b_{0}:=\inf\limits_{t\in[0,T]}b(t)$, we consider 
	\begin{equation}\label{maintranspde2}
		\left\{\begin{array}{l}
			\mathbb{D}_{(g)} \widehat{u}(t, \pi;a_{0},b_{0})_{m,k}+\pi_{m}^{2s}a_{0} \widehat{u}(t, \pi;a_{0},b_{0})_{m,k}+b_{0}\widehat{u}(t, \pi;a_{0},b_{0})_{m,k}=0, \\
			\widehat{u}(0, \pi;a_{0},b_{0})_{m,k}=\widehat{u}_{0}(\pi)_{m,k}, 
		\end{array}\right.
	\end{equation}
	that is,
	\begin{equation}\label{maintranspde2.1}
		\left\{\begin{array}{l}
			\mathbb{D}_{(g)} \widehat{u}(t, \pi;a_{0},b_{0})_{m,k}+\lambda_{\pi}^{s,m} \widehat{u}(t, \pi;a_{0},b_{0})_{m,k}=0,\quad t\in(0,T], \\
			\widehat{u}(0, \pi;a_{0},b_{0})_{m,k}=\widehat{u}_{0}(\pi)_{m,k},
		\end{array}\right.
	\end{equation}
	where $\lambda_{\pi}^{s,m}:=\pi_{m}^{2s}a_{0}+b_{0}>0$, since $\pi_{m}^{2s}>0$ and either $a_{0}> 0$ or $b_{0}> 0$.
	Using Remark \ref{3rdremark} for the Cauchy problem \eqref{maintranspde2.1}, we deduce that
	\begin{equation}\label{u0solution}
		\widehat{u}(t,\pi;a_{0},b_{0})_{m,k}=
		\widehat{u}_{0}(\pi)_{m,k}\widehat{u}_{\lambda_{\pi}^{s,m}}(t,\pi;a_{0},b_{0})_{m,k},\quad t\in(0,T],
	\end{equation}
	where $\widehat{u}_{\lambda_{\pi}^{s,m}}(t,\pi;a_{0},b_{0})_{m,k}$ solves the Cauchy problem \eqref{maintranspde2.1} with  $\widehat{u}(0,\pi;a_{0},b_{0})_{m,k}\equiv 1$.
	If we set 
	\begin{equation*}
		w(t, \pi;a,b)_{m,k}:=\widehat{u}(t,\pi;a_{0},b_{0})_{m,k}-\widehat{u}(t, \pi;a,b)_{m,k},\quad t\in[0,T],
	\end{equation*}
	then $w(0, \pi;a,b)_{m,k}\equiv 0$, and consequently, the equations \eqref{maintranspde1} and  \eqref{maintranspde2} give
	\begin{equation}\label{maintranspde3}
		\left\{\begin{array}{l}
			\mathbb{D}_{(g)} w(t, \pi;a,b)_{m,k}+\pi_{m}^{2s}a_{0} w(t, \pi;a,b)_{m,k}+b_{0}w(t, \pi;a,b)_{m,k}=\eta(t,\pi)_{m,k},~~t\in(0,T] \\
			w(0, \pi;a,b)_{m,k}\equiv0, 
		\end{array}\right.
	\end{equation}
	where
	\begin{equation*}
		\eta(t,\pi)_{m,k}:=\left[\pi_{m}^{2s}(a(t)-a_{0})+(b(t)-b_{0})\right]\widehat{u}(t, \pi;a,b)_{m,k},\quad t\in[0,T].
	\end{equation*}
	It is easy to note that
	\begin{equation}\label{signreg1}
		\eta(t,\pi)_{m,k} ~~\text{ holds the same sign as }~~ \widehat{u}(t, \pi;a,b)_{m,k},
	\end{equation}
	since $\pi_{m}^{2s}(a(t)-a_{0})+(b(t)-b_{0})\geq0$ for all $t\in[0,T]$. Combining \eqref{signreg} and \eqref{signreg1} along with
	Lemma \ref{regularity} for the Cauchy problem \eqref{maintranspde3}, we deduce that 
	\begin{equation}\label{signreg2}
		w(t, \pi;a,b)_{m,k}~~\text{ holds the same sign as }~~ \widehat{u}_{0}(\pi)_{m,k}.
	\end{equation}
	\sloppy  Now using the definition and regularity of $w(t, \pi;a,b)_{m,k}$, we conclude that
	\begin{equation*}
		\left\{\begin{array}{l}
			0 \leq \widehat{u}(t, \pi;a,b)_{m,k} \leq \widehat{u}(t,\pi;a_{0},b_{0})_{m,k}, \quad\text { if } \widehat{u}_{0}(\pi)_{m,k} \geq 0, \\
			\widehat{u}(t,\pi;a_{0},b_{0})_{m,k} \leq \widehat{u}(t, \pi;a,b)_{m,k} \leq 0,\quad \text { if }\widehat{u}_{0}(\pi)_{m,k} \leq 0.
		\end{array}\right.
	\end{equation*}
	Recalling Remark \ref{2ndremark}, i.e., $0\leq \widehat{u}_{\lambda_{\pi}^{s,m}}(t,\pi;a_{0},b_{0})_{m,k}\leq 1$,  for all $ t\in[0,T]$ and combining with the equation \eqref{u0solution}, we obtain
	\begin{equation}\label{u0est}
		|\widehat{u}(t, \pi;a,b)_{m,k}|	\leq|\widehat{u}(t, \pi;a_{0},b_{0})_{m,k}|
		=|\widehat{u}_{0}(\pi)_{m,k}\widehat{u}_{\lambda_{\pi}^{s,m}}(t,\pi;a_{0},b_{0})_{m,k}|
		\leq|\widehat{u}_{0}(\pi)_{m,k}|,
	\end{equation}
	for all $ t\in(0,T].$ If we multiply the above inequality by $\left(1+\pi_{m}^{2}\right)^{\frac{\gamma}{\nu}}\geq 1$ for some $\gamma\in\mathbb{R}$ and $\nu$ being the homogeneous degree of Rockland operator $\mathcal{R}$, then we get
	\begin{equation}\label{uu0est}
		\left|\left(1+\pi_{m}^{2}\right)^{\frac{\gamma}{\nu}}\widehat{u}(t, \pi;a,b)_{m,k}\right|	
		\leq \left|\left(1+\pi_{m}^{2}\right)^{\frac{\gamma}{\nu}}\widehat{u}_{0}(\pi)_{m,k}\right|, ~~\text{for all } t\in[0,T].
	\end{equation}
	Recalling the fact that for any Hilbert-Schmidt operator $S$, we have
	\begin{equation*}
		\|S\|^{2}_{\mathrm{HS}(\mathcal{H}_{\pi})}=\sum_{m,k}|(S\phi_{m},\phi_k)|^{2},
	\end{equation*}
	where $\{\phi_{1},\phi_{2},\dots\}$ is an orthonormal basis, we may compute the infinite sum of the inequalities \eqref{u0est}  and \eqref{uu0est} over $(m,k)\in\mathbb{N}^{2}$ to get the following estimate
	\begin{equation}\label{hsnorm}
		\|\widehat{u}(t, \pi;a,b)\|^{2}_{\mathrm{HS}(\mathcal{H}_{\pi})}\leq \|\widehat{u}_{0}(\pi)\|^{2}_{\mathrm{HS}(\mathcal{H}_{\pi})}
	\end{equation}
	and
	\begin{equation*}
		\left\|(1+\pi(\mathcal{R}))^{\frac{\gamma}{\nu}}\widehat{u}(t, \pi;a,b)\right\|^{2}_{\mathrm{HS}(\mathcal{H}_{\pi})}\leq \left\|(1+\pi(\mathcal{R}))^{\frac{\gamma}{\nu}}\widehat{u}_{0}(\pi)\right\|^{2}_{\mathrm{HS}(\mathcal{H}_{\pi})},
	\end{equation*}
	respectively.
	Now, by integrating both sides of above estimates against the Plancherel measure $\mu$ on $\widehat{\mathbb{G}}$ and using the Plancherel formula \eqref{planch} and the Sobolev norm \eqref{sobnorm}, we obtain
	\begin{equation}\label{est11}
		\left\|u(t,\cdot ; a,b)\right\|_{L^{2}(\mathbb{G})}
		\leq  \|u_{0}\|_{L^{2}(\mathbb{G})}~~\text{ and }~~  \left\|u(t,\cdot ; a,b)\right\|_{L^{2}_{\gamma}(\mathbb{G})}
		\leq  \|u_{0}\|_{L^{2}_{\gamma}(\mathbb{G})}.
	\end{equation}
	This proves the regularity estimate \eqref{MI01} and \eqref{MI02}. Similarily, using the equation \eqref{homeqnab}  with the above estimate, we get \begin{eqnarray*}
		\left\|\mathbb{D}_{(g)}u(t,\cdot ; a,b)\right\|_{L^{2}(\mathbb{G})}&=&\|-a(t)\mathcal{R}^{s}u(t,\cdot ; a,b)-b(t)u(t,\cdot ; a,b)\|_{L^{2}(\mathbb{G})}\nonumber\\
		&\leq&\|a\|_{C([0,T])}\|\mathcal{R}^{s}u(t,\cdot ; a,b)\|_{L^{2}(\mathbb{G})}+\|b\|_{C([0,T])}\|u(t,\cdot ; a,b)\|_{L^{2}(\mathbb{G})}\nonumber\\
		&\leq&C_{a,b}\left(\|\mathcal{R}^{s}u(t,\cdot ; a,b)\|_{L^{2}(\mathbb{G})}+\|u(t,\cdot ; a,b)\|_{L^{2}(\mathbb{G})}\right)\nonumber\\
		&\lesssim& \|(I+\mathcal{R})^{s}u(t,\cdot ; a,b)\|_{L^{2}(\mathbb{G})} \quad ( \text{using } \eqref{eqult})\nonumber\\
		&=&  \left\|u(t,\cdot ; a,b)\right\|_{L^{2}_{s\nu}(\mathbb{G})}\nonumber\\
		&\leq& \|u_{0}\|_{L^{2}_{s\nu}(\mathbb{G})} \quad (\text{using } \eqref{est11}),
	\end{eqnarray*}
	for all $t\in[0,T],$ where the constant $C_{a,b}:=\max\{\|a\|_{C([0,T])},\|b\|_{C([0,T])}\}$. This proves the regularity estimate \eqref{MI03}. Now consider 
	\begin{eqnarray*}
	&&\left\|\mathbb{D}_{(g)}u(t,\cdot ; a,b)\right\|^{2}_{L^{2}_{\gamma}(\mathbb{G})}\nonumber\\&=&\|(1+\mathcal{R})^{\frac{\gamma}{\nu}}\mathbb{D}_{(g)}u(t,\cdot ; a,b)\|^{2}_{L^{2}(\mathbb{G})}\nonumber\\
	&=& \|(1+\mathcal{R})^{\frac{\gamma}{\nu}}\left(-a(t)\mathcal{R}^{s}u(t,\cdot ; a,b)-b(t)u(t,\cdot ; a,b)\right)\|^{2}_{L^{2}(\mathbb{G})}\nonumber\\
	&\leq&C_{a,b}\int\limits_{\widehat{\mathbb{G}}}\left\|(1+\pi(\mathcal{R}))^{\frac{\gamma}{\nu}}\pi(\mathcal{R})^{s}\widehat{u}(t,\pi;a,b)+(1+\mathcal{R})^{\frac{\gamma}{\nu}}\widehat{u}(t,\pi ; a,b)\right\|^{2}_{\operatorname{HS}(\mathcal{H}_{\pi})}\mathrm{d}\mu(\pi)\nonumber\\
		&\lesssim&\int\limits_{\widehat{\mathbb{G}}}\left[\left\|(1+\pi(\mathcal{R}))^{\frac{\gamma+s\nu}{\nu}}\widehat{u}(t,\pi;a,b)\right\|^{2}_{\operatorname{HS}(\mathcal{H}_{\pi})}+\left\|(1+\mathcal{R})^{\frac{\gamma}{\nu}}\widehat{u}(t,\pi ; a,b)\right\|^{2}_{\operatorname{HS}(\mathcal{H}_{\pi})}\right]\mathrm{d}\mu(\pi)\nonumber\\&\leq&2\|u(t,\cdot ; a,b)\|^{2}_{L^{2}_{\gamma+s\nu}(\mathbb{G})}\quad (\text{using } \mathcal{R}\text{-Sobolev embeddings})\nonumber\\
		&\lesssim&\|u_{0}\|^{2}_{L^{2}_{\gamma+s\nu}(\mathbb{G})}\quad (\text{using } \eqref{est11}).
	\end{eqnarray*}
This concludes the regularity estimate \eqref{MI04} and hence the proof of Theorem \ref{mthm1}.
\end{proof}
Thus, we have established the regularity estimates for homogeneous Cauchy problem \eqref{mainpde}.

Let us now establish the regularity estimates for the non-homogeneous Cauchy problem \eqref{mainpdenonhom}.
In order to establish the regularity estimates stated in Theorem \ref{mthm2}, we divide the Cauchy problem \eqref{mainpdenonhom} into homogeneous and non-homogeneous parts.  

Let $u^{u_{0}}(t,x)$ and $u^{f}(t,x)$ represents the solutions of the following Cauchy problems:
\begin{equation}\label{homeq1}
	\left\{\begin{array}{l}
		\mathbb{D}_{(g)} u^{u_{0}}(t, x)+a(t)\mathcal{R}^{s} u^{u_{0}}(t, x)+b(t)u^{u_{0}}(t, x)=0, \quad (t,x) \in(0, T]\times \mathbb{G}, \\
		u^{u_{0}}(0, x)=u_{0}(x), \quad x \in  \mathbb{G},
	\end{array}\right.
\end{equation}
and
\begin{equation}\label{homeq2}
	\left\{\begin{array}{l}
	\mathbb{D}_{(g)} u^{f}(t, x)+a(t)\mathcal{R}^{s} u^{f}(t, x)+b(t)u^{f}(t, x)=f(t,x), \quad (t,x) \in(0, T]\times \mathbb{G}, \\
	u^{f}(0, x)=0, \quad x \in  \mathbb{G},
\end{array}\right.
\end{equation}
respectively. Now, we can write the solution $u(t,x)$ of the Cauchy problem \eqref{mainpdenonhom}  in terms of solutions of homogeneous Cauchy problem \eqref{homeq1} and non-homogeneous Cauchy problem \eqref{homeq2} as
\begin{equation}\label{eqnnn}
u(t,x)=u^{u_{0}}(t,x)+u^{f}(t,x),\quad (t,x)\in(0,T]\times\mathbb{G}.
\end{equation}
In order to determine the regularity estimate $u(t,x)$, we can apply Theorem \ref{mthm1} to obtain the estimates for the solution $u^{u_{0}}(t,x)$; however, we still require determining the regularity estimate for the solutions $u^{f}(t,x)$, which we will accomplish in the subsequent Lemma \ref{fnormm}.
\begin{lemma}\label{fnormm}
Let $\mathbb{G}$ be a graded Lie group and 
$\mathcal{R}$ be a positive Rockland operator of homogeneous degree $\nu$.  Assuming $0\leq a(t),b(t)\in C([0,T])$ such that  $\inf\limits_{t\in[0,T]}b(t)=b_{0}>0$, the solution $u^{f}(t,x)$ of the Cauchy problem \eqref{homeq2} satisfies the following estimates: 
\begin{enumerate}
	\item if  $f\in C([0,T];L^{2}(\mathbb{G}))$, then the solution $u^{f}$  satisfies the
	estimate  
	\begin{equation}\label{NMI01}
		\|u^{f}(t,\cdot)\|_{L^{2}(\mathbb{G})}\lesssim \|f\|_{C([0,T];L^{2}(\mathbb{G}))};
	\end{equation}
	\item if  $f\in C([0,T];L^{2}_{\gamma}(\mathbb{G}))$, for some $\gamma\in\mathbb{R}$, then the solution $u^{f}$  satisfies the
	estimate  
	\begin{equation}\label{NMI02}
		\|u^{f}(t,\cdot)\|_{L^{2}_{\gamma}(\mathbb{G})}\lesssim \|f\|_{C([0,T];L^{2}_{\gamma}(\mathbb{G}))};
	\end{equation}
	\item if  $f\in C([0,T];L^{2}_{s\nu}(\mathbb{G}))$, then the solution $u^{f}$  satisfies the
	estimate  
	\begin{equation}\label{NMI03}
		\|\mathbb{D}_{(g)}u^{f}(t,\cdot)\|_{L^{2}(\mathbb{G})}\lesssim \|f\|_{C([0,T];L^{2}_{s\nu}(\mathbb{G}))};
	\end{equation}
	\item if  $f\in C([0,T];L^{2}_{\gamma+s\nu}(\mathbb{G}))$, for some $\gamma\in\mathbb{R}$, then the solution $u^{f}$  satisfies the
	estimate  
	\begin{equation}\label{NMI04}
		\|\mathbb{D}_{(g)}u^{f}(t,\cdot)\|_{L^{2}_{\gamma}(\mathbb{G})}\lesssim \|f\|_{C([0,T];L^{2}_{\gamma+s\nu}(\mathbb{G}))};
	\end{equation}
\end{enumerate}
for all $t\in[0,T]$.
\end{lemma}
\begin{proof} Consider the non-homogeneous Cauchy problem
	\begin{equation*}
		\left\{\begin{array}{l}
			\mathbb{D}_{(g)} u^{f}(t, x;a,b)+a(t)\mathcal{R}^{s} u^{f}(t, x;a,b)+b(t)u^{f}(t, x;a,b)=f(t,x), \quad t \in(0, T], \\
			u^{f}(0, x;a,b)=0, \quad x \in  \mathbb{G}.
		\end{array}\right.
	\end{equation*}
Applying the group Fourier transform \eqref{grpft} using the relation \eqref{ftrock}, we get
	\begin{equation}\label{maintranspde11}
	\left\{\begin{array}{l}
	\mathbb{D}_{(g)} \widehat{u}^{f}(t, \pi;a,b)_{m,k}+\pi_{m}^{2s}a(t) \widehat{u}^{f}(t, \pi;a,b)_{m,k}+b(t) \widehat{u}^{f}(t, \pi;a,b)_{m,k}=\widehat{f}(t,\pi)_{m,k},   \\
		\widehat{u}^{f}(0, \pi;a,b)_{m,k}=0,
	\end{array}\right.
\end{equation}
for all $t\in(0,T]$. Similar to \eqref{maintranspde1}, we address each equation independently, by fixing $\pi\in\widehat{\mathbb{G}}$ and $(m,k)\in\mathbb{N}^{2}$. 	For $a_{0}:=\inf\limits_{t\in[0,T]}a(t)\geq 0$ and $b_{0}:=\inf\limits_{t\in[0,T]}b(t)>0$, we consider 
\begin{equation}\label{maintranspde22.}
	\left\{\begin{array}{l}
		\mathbb{D}_{(g)} \widehat{u}^{f}(t, \pi;a_{0},b_{0})_{m,k}+\pi_{m}^{2s}a_{0} \widehat{u}^{f}(t, \pi;a_{0},b_{0})_{m,k}+b_{0} \widehat{u}^{f}(t, \pi;a_{0},b_{0})_{m,k}=\widehat{f}(t,\pi)_{m,k},   \\
		\widehat{u}^{f}(0, \pi;a_{0},b_{0})_{m,k}=0,
	\end{array}\right.
\end{equation}
for all $t\in(0,T]$.
If we write 
\begin{equation*}
    \widehat{f}(t,\pi)_{m,k}=\mathrm{Re}     \widehat{f}(t,\pi)_{m,k}+i\mathrm{Im}    \widehat{f}(t,\pi)_{m,k},\quad \text{for all  }t\in[0,T],
\end{equation*}
 then we can denote 
\begin{equation}\label{fmaxresign}
    \widehat{f}^{\max}_{\mathrm{Re}}(t,\pi)_{m,k}:=\max\{0,\mathrm{Re}     \widehat{f}(t,\pi)_{m,k}\}=\frac{\mathrm{Re}     \widehat{f}(t,\pi)_{m,k}+|\mathrm{Re}     \widehat{f}(t,\pi)_{m,k}|}{2}\geq 0,     
\end{equation}
\begin{equation}\label{fminresign}
	\widehat{f}^{\min}_{\mathrm{Re}}(t,\pi)_{m,k}:=\min\{0,\mathrm{Re}     \widehat{f}(t,\pi)_{m,k}\}=\frac{\mathrm{Re}     \widehat{f}(t,\pi)_{m,k}-|\mathrm{Re}     \widehat{f}(t,\pi)_{m,k}|}{2}\leq 0,     
\end{equation}
\begin{equation}\label{fmaximsign}
	\widehat{f}^{\max}_{\mathrm{Im}}(t,\pi)_{m,k}:=\max\{0,\mathrm{Im}     \widehat{f}(t,\pi)_{m,k}\}=\frac{\mathrm{Im}     \widehat{f}(t,\pi)_{m,k}+|\mathrm{Im}     \widehat{f}(t,\pi)_{m,k}|}{2}\geq 0,     
\end{equation}
and
\begin{equation}\label{fminimsign}
	\widehat{f}^{\min}_{\mathrm{Im}}(t,\pi)_{m,k}:=\min\{0,\mathrm{Im}     \widehat{f}(t,\pi)_{m,k}\}=\frac{\mathrm{Im}     \widehat{f}(t,\pi)_{m,k}-|\mathrm{Im}     \widehat{f}(t,\pi)_{m,k}|}{2}\leq 0.  
\end{equation}
 In light of the situation that $\mathrm{Re}\widehat{f}=\widehat{f}^{\max}_{\mathrm{Re}}+\widehat{f}^{\min}_{\mathrm{Re}}$ and $\mathrm{Im}\widehat{f}=\widehat{f}^{\max}_{\mathrm{Im}}+\widehat{f}^{\min}_{\mathrm{Im}}$, the Cauchy problems \eqref{maintranspde11} and \eqref{maintranspde22.} can be divided into the subsequent Cauchy problems:
	\begin{equation}\label{fmaxa:Re}
	\left\{\begin{array}{l}
		\mathbb{D}_{(g)} \widehat{u}^{f,\max}_{\mathrm{Re}}(t, \pi;a,b)_{m,k}+[\pi_{m}^{2s}a(t)+b(t)] \widehat{u}^{f,\max}_{\mathrm{Re}}(t, \pi;a,b)_{m,k}=\widehat{f}^{\max}_{\mathrm{Re}}(t,\pi)_{m,k},   \\
		\widehat{u}^{f,\max}_{\mathrm{Re}}(0, \pi;a,b)_{m,k}=0, 
	\end{array}\right.
\end{equation}
\begin{equation}\label{fmaxa0:Re}
	\left\{\begin{array}{l}
		\mathbb{D}_{(g)} \widehat{u}^{f,\max}_{\mathrm{Re}}(t, \pi;a_{0},b_{0})_{m,k}+[\pi_{m}^{2s}a_{0}+b_{0}] \widehat{u}^{f,\max}_{\mathrm{Re}}(t, \pi;a_{0},b_{0})_{m,k}=\widehat{f}^{\max}_{\mathrm{Re}}(t,\pi)_{m,k},   \\
		\widehat{u}^{f,\max}_{\mathrm{Re}}(0, \pi;a_{0},b_{0})_{m,k}=0, 
	\end{array}\right.
\end{equation}
\begin{equation}\label{fmina:Re}
	\left\{\begin{array}{l}
		\mathbb{D}_{(g)} \widehat{u}^{f,\min}_{\mathrm{Re}}(t, \pi;a,b)_{m,k}+[\pi_{m}^{2s}a(t)+b(t)] \widehat{u}^{f,\min}_{\mathrm{Re}}(t, \pi;a,b)_{m,k}=\widehat{f}^{\min}_{\mathrm{Re}}(t,\pi)_{m,k},   \\
		\widehat{u}^{f,\min}_{\mathrm{Re}}(0, \pi;a,b)_{m,k}=0, 
	\end{array}\right.
\end{equation}
\begin{equation}\label{fmina0:Re}
	\left\{\begin{array}{l}
		\mathbb{D}_{(g)} \widehat{u}^{f,\min}_{\mathrm{Re}}(t, \pi;a_{0},b_{0})_{m,k}+[\pi_{m}^{2s}a_{0}+b_{0}] \widehat{u}^{f,\min}_{\mathrm{Re}}(t, \pi;a_{0},b_{0})_{m,k}=\widehat{f}^{\min}_{\mathrm{Re}}(t,\pi)_{m,k},   \\
		\widehat{u}^{f,\min}_{\mathrm{Re}}(0, \pi;a_{0},b_{0})_{m,k}=0,
	\end{array}\right.
\end{equation}
	\begin{equation}\label{fmaxa:Im}
	\left\{\begin{array}{l}
		\mathbb{D}_{(g)} \widehat{u}^{f,\max}_{\mathrm{Im}}(t, \pi;a,b)_{m,k}+[\pi_{m}^{2s}a(t)+b(t)] \widehat{u}^{f,\max}_{\mathrm{Im}}(t, \pi;a,b)_{m,k}=\widehat{f}^{\max}_{\mathrm{Im}}(t,\pi)_{m,k},   \\
		\widehat{u}^{f,\max}_{\mathrm{Im}}(0, \pi;a,b)_{m,k}=0, 
	\end{array}\right.
\end{equation}
\begin{equation}\label{fmaxa0:Im}
	\left\{\begin{array}{l}
		\mathbb{D}_{(g)} \widehat{u}^{f,\max}_{\mathrm{Im}}(t, \pi;a_{0},b_{0})_{m,k}+[\pi_{m}^{2s}a_{0}+b_{0}] \widehat{u}^{f,\max}_{\mathrm{Im}}(t, \pi;a_{0},b_{0})_{m,k}=\widehat{f}^{\max}_{\mathrm{Im}}(t,\pi)_{m,k},   \\
		\widehat{u}^{f,\max}_{\mathrm{Im}}(0, \pi;a_{0},b_{0})_{m,k}=0, 
	\end{array}\right.
\end{equation}
\begin{equation}\label{fmina:Im}
	\left\{\begin{array}{l}
		\mathbb{D}_{(g)} \widehat{u}^{f,\min}_{\mathrm{Im}}(t, \pi;a,b)_{m,k}+[\pi_{m}^{2s}a(t)+b(t)] \widehat{u}^{f,\min}_{\mathrm{Im}}(t, \pi;a,b)_{m,k}=\widehat{f}^{\min}_{\mathrm{Im}}(t,\pi)_{m,k},   \\
		\widehat{u}^{f,\min}_{\mathrm{Im}}(0, \pi;a,b)_{m,k}=0, 
	\end{array}\right.
\end{equation}
and
\begin{equation}\label{fmina0:Im}
	\left\{\begin{array}{l}
		\mathbb{D}_{(g)} \widehat{u}^{f,\min}_{\mathrm{Im}}(t, \pi;a_{0},b_{0})_{m,k}+[\pi_{m}^{2s}a_{0}+b_{0}] \widehat{u}^{f,\min}_{\mathrm{Im}}(t, \pi;a_{0},b_{0})_{m,k}=\widehat{f}^{\min}_{\mathrm{Im}}(t,\pi)_{m,k},   \\
		\widehat{u}^{f,\min}_{\mathrm{Im}}(0, \pi;a_{0},b_{0})_{m,k}=0, 
	\end{array}\right.
\end{equation}
for all $t\in[0,T]$. By applying Lemma \ref{regularity} to the Cauchy problems \eqref{fmaxa:Re}, \eqref{fmina:Re}, \eqref{fmaxa:Im},   and \eqref{fmina:Im}, we can now conclude that
\begin{equation}\label{cond1w}
	\left\{\begin{array}{l}
		\widehat{u}^{f, \max }_{\mathrm{Re}}(t,\pi ; a,b)_{m,k}\geq 0, \quad \text{since}~~ \widehat{f}^{\max}_{\mathrm{Re}}(t,\pi)_{m,k}\geq 0;\\
		\widehat{u}^{f, \min }_{\mathrm{Re}}(t,\pi ; a,b)_{m,k}\leq 0, \quad \text{since}~~ \widehat{f}^{\min}_{\mathrm{Re}}(t,\pi)_{m,k}\leq 0;\\
		\widehat{u}^{f, \max }_{\mathrm{Im}}(t,\pi ; a,b)_{m,k}\geq0, \quad \text{since}~~\widehat{f}^{\max}_{\mathrm{Im}}(t,\pi)_{m,k}\geq 0; ~\text{and}\\
		\widehat{u}^{f, \min}_{\mathrm{Im}}(t,\pi ; a,b)_{m,k}\leq 0 ,\quad \text{since}~~ \widehat{f}^{\min}_{\mathrm{Im}}(t,\pi)_{m,k}\leq 0;
	\end{array}\right.
\end{equation}
for all $t\in[0,T]$. Further, if we denote
\begin{equation*}
w^{f,\max}_{\mathrm{Re}}(t, \pi;a,b)_{m,k}:=\widehat{u}^{f,\max}_{\mathrm{Re}}(t,\pi;a_{0},b_{0})_{m,k}-\widehat{u}^{f,\max}_{\mathrm{Re}}(t, \pi;a,b)_{m,k},
\end{equation*}
\begin{equation*}
w^{f,\min}_{\mathrm{Re}}(t, \pi;a,b)_{m,k}:=\widehat{u}^{f,\min}_{\mathrm{Re}}(t,\pi;a_{0},b_{0})_{m,k}-\widehat{u}^{f,\min}_{\mathrm{Re}}(t, \pi;a,b)_{m,k},
\end{equation*}
\begin{equation*}
	w^{f,\max}_{\mathrm{Im}}(t, \pi;a,b)_{m,k}:=\widehat{u}^{f,\max}_{\mathrm{Im}}(t,\pi;a_{0},b_{0})_{m,k}-\widehat{u}^{f,\max}_{\mathrm{Im}}(t, \pi;a,b)_{m,k},
\end{equation*}
and
\begin{equation*}
	w^{f,\min}_{\mathrm{Im}}(t, \pi;a,b)_{m,k}:=\widehat{u}^{f,\min}_{\mathrm{Im}}(t,\pi;a_{0},b_{0})_{m,k}-\widehat{u}^{f,\min}_{\mathrm{Im}}(t, \pi;a,b)_{m,k},
\end{equation*}
then using the equations \eqref{fmaxa:Re}-\eqref{fmina0:Im}, we deduce that
\begin{equation}\label{wmax:Re}
	\left\{\begin{array}{l}
\mathbb{D}_{(g)}w^{f,\max}_{\mathrm{Re}}(t, \pi;a,b)_{m,k}+[\pi_{m}^{2s}a_{0}+b_{0}] w^{f,\max}_{\mathrm{Re}}(t, \pi;a,b)_{m,k}=\chi_{\mathrm{Re}}^{\max}(t,\pi)_{m,k},~~t\in(0,T],\\
		w^{f,\max}_{\mathrm{Re}}(t, \pi;a,b)_{m,k}=0, 
	\end{array}\right.
\end{equation}
where 
\begin{equation}\label{chi:max:Re}
	\chi_{\mathrm{Re}}^{\max}(t,\pi)_{m,k}:=\left[\pi_{m}^{2s}(a(t)-a_{0})+(b(t)-b_{0})\right]\widehat{u}^{f,\max}_{\mathrm{Re}}(t, \pi;a,b)_{m,k},
\end{equation}
\begin{equation}\label{wmin:Re}
	\left\{\begin{array}{l}
		\mathbb{D}_{(g)}w^{f,\min}_{\mathrm{Re}}(t, \pi;a,b)_{m,k}+[\pi_{m}^{2s}a_{0}+b_{0}] w^{f,\min}_{\mathrm{Re}}(t, \pi;a,b)_{m,k}=\chi_{\mathrm{Re}}^{\min}(t,\pi)_{m,k},~~t\in(0,T],\\
		w^{f,\min}_{\mathrm{Re}}(t, \pi;a,b)_{m,k}=0, 
	\end{array}\right.
\end{equation}
where 
\begin{equation}\label{chi:min:Re}
	\chi_{\mathrm{Re}}^{\min}(t,\pi)_{m,k}:=\left[\pi_{m}^{2s}(a(t)-a_{0})+(b(t)-b_{0})\right]\widehat{u}^{f,\min}_{\mathrm{Re}}(t, \pi;a,b)_{m,k},
\end{equation}
\begin{equation}\label{wmax:Im}
	\left\{\begin{array}{l}
		\mathbb{D}_{(g)}w^{f,\max}_{\mathrm{Im}}(t, \pi;a,b)_{m,k}+[\pi_{m}^{2s}a_{0}+b_{0}] w^{f,\max}_{\mathrm{Im}}(t, \pi;a,b)_{m,k}=\chi_{\mathrm{Im}}^{\max}(t,\pi)_{m,k},~~t\in(0,T],\\
		w^{f,\max}_{\mathrm{Im}}(t, \pi;a,b)_{m,k}=0, 
	\end{array}\right.
\end{equation}
where 
\begin{equation}\label{chi:max:Im}
	\chi_{\mathrm{Im}}^{\max}(t,\pi)_{m,k}:=\left[\pi_{m}^{2s}(a(t)-a_{0})+(b(t)-b_{0})\right]\widehat{u}^{f,\max}_{\mathrm{Im}}(t, \pi;a,b)_{m,k},\quad \text{and},
\end{equation}
\begin{equation}\label{wmin:Im}
	\left\{\begin{array}{l}
		\mathbb{D}_{(g)}w^{f,\min}_{\mathrm{Im}}(t, \pi;a,b)_{m,k}+[\pi_{m}^{2s}a_{0}+b_{0}] w^{f,\min}_{\mathrm{Im}}(t, \pi;a,b)_{m,k}=\chi_{\mathrm{Im}}^{\min}(t,\pi)_{m,k},~~t\in(0,T],\\
		w^{f,\min}_{\mathrm{Im}}(t, \pi;a,b)_{m,k}=0, 
	\end{array}\right.
\end{equation}
where 
\begin{equation}\label{chi:min:Im}
	\chi_{\mathrm{Im}}^{\min}(t,\pi)_{m,k}:=\left[\pi_{m}^{2s}(a(t)-a_{0})+(b(t)-b_{0})\right]\widehat{u}^{f,\max}_{\mathrm{Im}}(t, \pi;a,b)_{m,k}.
\end{equation}
\sloppy Combining the fact that $\pi_{m}^{2s}(a(t)-a_{0})+(b(t)-b_{0})\geq0$ for all $t\in[0,T]$ with the equations \eqref{chi:max:Re}, \eqref{chi:min:Re}, \eqref{chi:max:Im}, and \eqref{chi:min:Im},  it is easy to conclude that 
\begin{equation}\label{cond1eta}
	\left\{\begin{array}{l}
	\chi_{\mathrm{Re}}^{\max}(t,\pi)_{m,k}\geq	0,\quad \text{since}~~	\widehat{u}^{f, \max }_{\mathrm{Re}}(t,\pi ; a,b)_{m,k}\geq 0,  \\
	\chi_{\mathrm{Re}}^{\min}(t,\pi)_{m,k}\leq 0,	\quad \text{since}~~	\widehat{u}^{f, \min }_{\mathrm{Re}}(t,\pi ; a,b)_{m,k}\leq 0, \\
	\chi_{\mathrm{Im}}^{\max}(t,\pi)_{m,k}\geq 0,	\quad \text{since}~~	\widehat{u}^{f, \max }_{\mathrm{Im}}(t,\pi ; a,b)_{m,k}\geq0, \quad ~\text{and}\\
	\chi_{\mathrm{Im}}^{\min}(t,\pi)_{m,k}	\leq 0,\quad \text{since}~~	\widehat{u}^{f, \min}_{\mathrm{Im}}(t,\pi ; a,b)_{m,k}\leq 0 .
	\end{array}\right.
\end{equation}
Further, utilizing the Lemma \ref{regularity} along with the relations \eqref{cond1w} and \eqref{cond1eta} for the Cauchy problems \eqref{wmax:Re}, \eqref{wmin:Re}, \eqref{wmax:Im}, and  \eqref{wmin:Im}, we extract the following conclusion for the solutions
\begin{equation}\label{cond1sol}
	\left\{\begin{array}{l}
	w^{f,\max}_{\mathrm{Re}}(t, \pi;a,b)_{m,k}\geq 0,\quad \text{since}~~		\chi_{\mathrm{Re}}^{\max}(t,\pi)_{m,k}\geq 0,  \\
	w^{f,\min}_{\mathrm{Re}}(t, \pi;a,b)_{m,k}\leq 0,	\quad \text{since}~~	\chi_{\mathrm{Re}}^{\min}(t,\pi)_{m,k}\leq 0, \\
		w^{f,\max}_{\mathrm{Im}}(t, \pi;a,b)_{m,k}\geq 0,	\quad \text{since}~~	\chi_{\mathrm{Im}}^{\max}(t,\pi)_{m,k}\geq 0, \quad ~\text{and}\\
		w^{f,\min}_{\mathrm{Im}}(t, \pi;a,b)_{m,k}	\leq 0,\quad \text{since}~~	\chi_{\mathrm{Im}}^{\min}(t,\pi)_{m,k}\leq 0,
	\end{array}\right.
\end{equation}
for all $t\in[0,T]$. Now using the definition and regularity of $w^{f,\max}_{\operatorname{Re}},w^{f,\min}_{\operatorname{Re}},w^{f,\max}_{\operatorname{Im}},$ and $w^{f,\min}_{\operatorname{Im}}$, we conclude that
\begin{equation}\label{cond1}
\left\{\begin{array}{l}
0 \leq \widehat{u}^{f, \max }_{\mathrm{Re}}(t,\pi ; a,b)_{m,k} \leq \widehat{u}^{f, \max }_{\mathrm{Re}}\left(t ,\pi; a_{0},b_{0}\right)_{m,k},\quad \text { on }[0, T], \\
\widehat{u}^{f, \min }_{\mathrm{Re}}\left(t,\pi ; a_{0},b_{0}\right)_{m,k} \leq \widehat{u}^{f, \min }_{\mathrm{Re}}(t ,\pi; a,b)_{m,k} \leq 0,\quad \text { on }[0, T],\\
0 \leq \widehat{u}^{f, \max }_{\mathrm{Im}}(t,\pi ; a,b)_{m,k} \leq \widehat{u}^{f, \max }_{\mathrm{Im}}\left(t ,\pi; a_{0},b_{0}\right)_{m,k},\quad \text { on }[0, T],~~\text{and} \\
\widehat{u}^{f, \min }_{\mathrm{Im}}\left(t,\pi ; a_{0},b_{0}\right)_{m,k} \leq \widehat{u}^{f, \min }_{\mathrm{Im}}(t ,\pi; a,b)_{m,k} \leq 0,\quad \text { on }[0, T].
\end{array}\right.
\end{equation}
The moment Remark \ref{3rdremark} is applied to the Cauchy problems \eqref{maintranspde22.}, \eqref{fmaxa0:Re}, \eqref{fmina0:Re}, \eqref{fmaxa0:Im}, and \eqref{fmina0:Im}, the resulting solutions are as follows:
\begin{equation*}
	\widehat{u}^{f}(t,\pi;a_{0},b_{0})_{m,k}=-\frac{1}{\lambda_{\pi}^{s,m}} \int_{0}^{t} \widehat{f}(s,\pi)_{m,k}\frac{\partial}{\partial s}w_{\lambda_{\pi}^{s,m}}(t-s,\pi)_{m,k}  \mathrm{d} s,
	\end{equation*}
 \begin{equation*}
	\widehat{u}^{f,\max}_{\mathrm{Re}}(t,\pi;a_{0},b_{0})_{m,k}=-\frac{1}{\lambda_{\pi}^{s,m}} \int_0^t \widehat{f}^{\max}_{\mathrm{Re}}(s,\pi)_{m,k}\frac{\partial}{\partial s}w_{\lambda_{\pi}^{s,m}}(t-s,\pi)_{m,k}  \mathrm{d} s,
	\end{equation*}
	 \begin{equation*}
		\widehat{u}^{f,\min}_{\mathrm{Re}}(t,\pi;a_{0},b_{0})_{m,k}=-\frac{1}{\lambda_{\pi}^{s,m}} \int_0^t \widehat{f}^{\min}_{\mathrm{Re}}(s,\pi)_{m,k}\frac{\partial}{\partial s}w_{\lambda_{\pi}^{s,m}}(t-s,\pi)_{m,k}  \mathrm{d} s,
	\end{equation*}
	 \begin{equation*}
		\widehat{u}^{f,\max}_{\mathrm{Im}}(t,\pi;a_{0},b_{0})_{m,k}=-\frac{1}{\lambda_{\pi}^{s,m}} \int_0^t \widehat{f}^{\max}_{\mathrm{Im}}(s,\pi)_{m,k}\frac{\partial}{\partial s}w_{\lambda_{\pi}^{s,m}}(t-s,\pi)_{m,k}  \mathrm{d} s,
	\end{equation*}
 and
 \begin{equation*}
	\widehat{u}^{f,\min}_{\mathrm{Im}}(t,\pi;a_{0},b_{0})_{m,k}=-\frac{1}{\lambda_{\pi}^{s,m}} \int_0^t \widehat{f}^{\min}_{\mathrm{Im}}(s,\pi)_{m,k}\frac{\partial}{\partial s}w_{\lambda_{\pi}^{s,m}}(t-s,\pi)_{m,k}  \mathrm{d} s,
	\end{equation*}
 respectively, where $\lambda_{\pi}^{s,m}:=\pi_{m}^{2s}a_{0}+b_{0}$.
Putting together the equalities
\begin{eqnarray*}
\widehat{u}^{f}(t, \pi;a,b)_{m,k}&=&\operatorname{Re}\widehat{u}^{f}(t, \pi;a,b)_{m,k}+i\operatorname{Im}\widehat{u}^{f}(t, \pi;a,b)_{m,k}\\
&=&\widehat{u}^{f,\max}_{\mathrm{Re}}(t, \pi;a,b)_{m,k}+\widehat{u}^{f,\min}_{\mathrm{Re}}(t, \pi;a,b)_{m,k}\\&+&i\left(\widehat{u}^{f,\max}_{\mathrm{Im}}(t, \pi;a,b)_{m,k}+\widehat{u}^{f,\min}_{\mathrm{Im}}(t, \pi;a,b)_{m,k}\right),
\end{eqnarray*}
with the above equations, relations \eqref{fmaxresign}-\eqref{fminimsign}, inequality \eqref{cond1} and the hypothesis $f\in C([0,T];L^{2}(\mathbb{G}))$, we get
\begin{eqnarray}\label{festi}
    &&|\widehat{u}^{f}(t, \pi;a,b)_{m,k}|^{2}\nonumber\\&\leq&|\widehat{u}^{f,\max}_{\mathrm{Re}}(t, \pi;a,b)_{m,k}+\widehat{u}^{f,\min}_{\mathrm{Re}}(t, \pi;a,b)_{m,k}|^{2}\nonumber\\&&+|\widehat{u}^{f,\max}_{\mathrm{Im}}(t, \pi;a,b)_{m,k}+\widehat{u}^{f,\min}_{\mathrm{Im}}(t, \pi;a,b)_{m,k}|^{2}\nonumber\\
     &\leq&\left[\frac{1}{\lambda_{\pi}^{s,m}} \int_0^t \left( |\widehat{f}^{\max}_{\mathrm{Re}}(s,\pi)_{m,k}|+|\widehat{f}^{\min}_{\mathrm{Re}}(s,\pi)_{m,k}|\right)\left|\frac{\partial}{\partial s}w_{\lambda_{\pi}^{s,m}}(t-s,\pi)_{m,k} \right| \mathrm{d} s\right]^{2}\nonumber\\
    &&+\left[\frac{1}{\lambda_{\pi}^{s,m}} \int_0^t \left( |\widehat{f}^{\max}_{\mathrm{Im}}(s,\pi)_{m,k}|+|\widehat{f}^{\min}_{\mathrm{Im}}(s,\pi)_{m,k}|\right)\left|\frac{\partial}{\partial s}w_{\lambda_{\pi}^{s,m}}(t-s,\pi)_{m,k} \right| \mathrm{d} s\right]^{2}\nonumber\\
    &=&\left[\frac{1}{\lambda_{\pi}^{s,m}} \int_0^t \left( \widehat{f}^{\max}_{\mathrm{Re}}(s,\pi)_{m,k}-\widehat{f}^{\min}_{\mathrm{Re}}(s,\pi)_{m,k}\right)\left|\frac{\partial}{\partial s}w_{\lambda_{\pi}^{s,m}}(t-s,\pi)_{m,k} \right| \mathrm{d} s\right]^{2}\nonumber\\
    &&+\left[\frac{1}{\lambda_{\pi}^{s,m}} \int_0^t \left( \widehat{f}^{\max}_{\mathrm{Im}}(s,\pi)_{m,k}-\widehat{f}^{\min}_{\mathrm{Im}}(s,\pi)_{m,k}\right)\left|\frac{\partial}{\partial s}w_{\lambda_{\pi}^{s,m}}(t-s,\pi)_{m,k} \right| \mathrm{d} s\right]^{2}\nonumber\\
    &=&\left[\frac{1}{\lambda_{\pi}^{s,m}} \int_0^t  \left|\mathrm{Re}\widehat{f}(s,\pi)_{m,k}\right|\left|\frac{\partial}{\partial s}w_{\lambda_{\pi}^{s,m}}(t-s,\pi)_{m,k} \right| \mathrm{d} s\right]^{2}\nonumber\\
    &&+\left[\frac{1}{\lambda_{\pi}^{s,m}} \int_0^t \left|\mathrm{Im}\widehat{f}(s,\pi)_{m,k}\right|\left|\frac{\partial}{\partial s}w_{\lambda_{\pi}^{s,m}}(t-s,\pi)_{m,k} \right| \mathrm{d} s\right]^{2}\nonumber\\
    &\leq&\left[\frac{1}{\lambda_{\pi}^{s,m}} \sup_{s\in[0,T]}\left\{\left|\mathrm{Re}\widehat{f}(s,\pi)_{m,k}\right|\right\}\int_0^t  \left|\frac{\partial}{\partial s}w_{\lambda_{\pi}^{s,m}}(t-s,\pi)_{m,k} \right| \mathrm{d} s\right]^{2}\nonumber\\
    &&+\left[\frac{1}{\lambda_{\pi}^{s,m}} \sup_{s\in[0,T]}\left\{\left|\mathrm{Im}\widehat{f}(s,\pi)_{m,k}\right|\right\}\int_0^t \left|\frac{\partial}{\partial s}w_{\lambda_{\pi}^{s,m}}(t-s,\pi)_{m,k} \right| \mathrm{d} s\right]^{2}\nonumber\\
    &=&\left[\frac{1}{\lambda_{\pi}^{s,m}} \sup_{s\in[0,T]}\left\{\left|\mathrm{Re}\widehat{f}(s,\pi)_{m,k}\right|\right\}\left(1-w_{\lambda_{\pi}^{s,m}}(t,\pi)_{m,k}\right)\right]^{2}\quad (\because w^{\prime}_{\lambda_{\pi}^{s,m}}\leq 0)\nonumber\\
    &&+\left[\frac{1}{\lambda_{\pi}^{s,m}} \sup_{s\in[0,T]}\left\{\left|\mathrm{Im}\widehat{f}(s,\pi)_{m,k}\right|\right\}\left(1-w_{\lambda_{\pi}^{s,m}}(t,\pi)_{m,k}\right)\right]^{2}\quad (\because w_{\lambda_{\pi}^{s,m}}(0,\pi)=1)\nonumber\\
     &\leq&\frac{1}{(\lambda_{\pi}^{s,m})^{2}}\sup_{s\in[0,T]}\left\{\left|\mathrm{Re}\widehat{f}(s,\pi)_{m,k}\right|^{2}+\left|\mathrm{Im}\widehat{f}(s,\pi)_{m,k}\right|^{2}\right\}\quad (\because a\leq w_{\lambda_{\pi}^{s,m}}(t,\pi)_{m,k}\leq1)\nonumber\\
     &=&\frac{1}{(\lambda_{\pi}^{s,m})^{2}}\sup_{s\in[0,T]}\left\{\left|\widehat{f}(s,\pi)_{m,k}\right|^{2}\right\}\nonumber\\
      &\leq& \frac{\|\widehat{f}(\cdot,\pi)_{m,k}\|^{2}_{C([0,T])}}{b_{0}^2},
\end{eqnarray}
for all $t\in[0,T]$, where the last inequality follows from the fact that $\lambda_{\pi}^{s,m}=\pi_{m}^{2s}a_{0}+b_{0}>b_{0}>0$.  If we multiply the estimate \eqref{festi} by $\left(1+\pi_{m}^{2}\right)^{\frac{2\gamma}{\nu}}$ for some $\gamma\in\mathbb{R}$ and $\nu$ being the homogeneous degree of Rockland operator $\mathcal{R}$, then we get
\begin{equation}\label{festim}
\left|\left(1+\pi_{m}^{2}\right)^{\frac{\gamma}{\nu}}\widehat{u}^{f}(t, \pi;a,b)_{m,k}\right|^{2}\lesssim \left\|\left(1+\pi_{m}^{2}\right)^{\frac{2\gamma}{\nu}}\widehat{f}(\cdot,\pi)_{m,k}\right\|^{2}_{C([0,T])}.
\end{equation}
 Similar to the estimate \eqref{hsnorm}, the estimate \eqref{festi} and 
 \eqref{festim} gives
\begin{equation*}
	\left\|\widehat{u}^{f}(t,\pi ; a,b)\right\|^{2}_{\operatorname{HS}(\mathcal{H}_{\pi})}
\lesssim \|\widehat{f}(\cdot,\pi)\|^{2}_{C([0,T];\operatorname{HS}(\mathcal{H}_{\pi}))},
\end{equation*}
and
\begin{equation*}
\left\|\left(1+\pi(\mathcal{R})\right)^{\frac{\gamma}{\nu}}\widehat{u}^{f}(t,\pi ; a,b)\right\|^{2}_{\operatorname{HS}(\mathcal{H}_{\pi})}
\lesssim \left\|\left(1+\pi(\mathcal{R})\right)^{\frac{\gamma}{\nu}}\widehat{f}(\cdot,\pi)\right\|^{2}_{C([0,T];\operatorname{HS}(\mathcal{H}_{\pi}))}.
\end{equation*}
 Now, by integrating both sides of above estimates against the Plancherel measure $\mu$ on $\widehat{\mathbb{G}}$ and using the Plancherel formula \eqref{planch} and the Sobolev norm \eqref{sobnorm}, we obtain
  \begin{equation}\label{eqnn2}
 	\left\|u^{f}(t,\cdot ; a,b)\right\|_{L^{2}(\mathbb{G})}
 	\leq  \|f\|_{C([0,T];L^{2}(\mathbb{G}))}~~\text{ and }~~  \left\|u^{f}(t,\cdot ; a,b)\right\|_{L^{2}_{\gamma}(\mathbb{G})}
 	\leq  \|f\|_{C([0,T];L^{2}_{\gamma}(\mathbb{G}))}.
 \end{equation}
This proves the regularity estimates \eqref{NMI01} and \eqref{NMI02}.  Similarily, using the equation \eqref{homeq2}  with above estimate, we get
\begin{eqnarray*}
	&&\left\|\mathbb{D}_{(g)}u^{f}(t,\cdot ; a,b)\right\|_{L^{2}(\mathbb{G})}\nonumber\\&=&\|-a(t)\mathcal{R}^{s}u^{f}(t,\cdot ; a,b)-b(t)u^{f}(t,\cdot ; a,b)+f(t,\cdot)\|_{L^{2}(\mathbb{G})}\nonumber\\
	&\leq&\|a\|_{C([0,T])}\|\mathcal{R}^{s}u^{f}(t,\cdot ; a,b)\|_{L^{2}(\mathbb{G})}+\|b\|_{C([0,T])}\|u^{f}(t,\cdot ; a,b)\|_{L^{2}(\mathbb{G})}+\|f(t,\cdot)\|_{L^{2}(\mathbb{G})}\nonumber\\
	&\leq&C_{a,b}\left(\|u^{f}(t,\cdot ; a,b)\|_{L^{2}_{s\nu}(\mathbb{G})}+\|u^{f}(t,\cdot ; a,b)\|_{L^{2}(\mathbb{G})}+\|f(t,\cdot)\|_{L^{2}(\mathbb{G})}\right)\nonumber\\
	&\lesssim&\|u^{f}(t,\cdot ; a,b)\|_{L^{2}_{s\nu}(\mathbb{G})}+\|f(t,\cdot)\|_{L^{2}(\mathbb{G})}\quad (\text{using } \eqref{eqult}  )\nonumber\\
&\leq&	 \|f\|_{C([0,T];L^{2}_{s\nu}(\mathbb{G}))}+\|f\|_{C([0,T];L^{2}(\mathbb{G}))}\quad (\text{using } \eqref{eqnn2})\nonumber\\
&\lesssim& \|f\|_{C([0,T];L^{2}_{s\nu}(\mathbb{G}))} \quad (\text{using } \eqref{eqult}  ), 
\end{eqnarray*}
for all $t\in[0,T],$ with the constant $C_{a,b}:=\max\{\|a\|_{C([0,T])},\|b\|_{C([0,T])},1\}$. This proves the regularity estimate \eqref{NMI03}.
Now consider 
\begin{eqnarray*}
	&&\left\|\mathbb{D}_{(g)}u^{f}(t,\cdot ; a,b)\right\|^{2}_{L^{2}_{\gamma}(\mathbb{G})}\nonumber\\&=&\left\|(1+\mathcal{R})^{\frac{\gamma}{\nu}}\mathbb{D}_{(g)}u^{f}(t,\cdot ; a,b)\right\|^{2}_{L^{2}(\mathbb{G})}\nonumber\\
	&=& \left\|(1+\mathcal{R})^{\frac{\gamma}{\nu}}\left(-a(t)\mathcal{R}^{s}u^{f}(t,\cdot ; a,b)-b(t)u^{f}(t,\cdot ; a,b)+f(t,\cdot)\right)\right\|^{2}_{L^{2}(\mathbb{G})}\nonumber\\
	&\leq&C_{a,b}\int\limits_{\widehat{\mathbb{G}}}\left\|(1+\pi(\mathcal{R}))^{\frac{\gamma}{\nu}}\left(\pi(\mathcal{R})^{s}\widehat{u}^{f}(t,\pi;a,b)+\widehat{u}^{f}(t,\pi ; a,b)+\widehat{f}(t,\cdot)\right)\right\|^{2}_{\operatorname{HS}(\mathcal{H}_{\pi})}\mathrm{d}\mu(\pi)\nonumber\\
	&\lesssim&\int\limits_{\widehat{\mathbb{G}}}\left[\left\|(1+\pi(\mathcal{R}))^{\frac{\gamma+s\nu}{\nu}}\widehat{u}^{f}(t,\pi;a,b)\right\|^{2}_{\operatorname{HS}(\mathcal{H}_{\pi})}+\left\|(1+\mathcal{R})^{\frac{\gamma}{\nu}}\widehat{u}^{f}(t,\pi ; a,b)\right\|^{2}_{\operatorname{HS}(\mathcal{H}_{\pi})}\right]\mathrm{d}\mu(\pi)\nonumber\\
	&+&\int\limits_{\widehat{\mathbb{G}}}\left\|(1+\mathcal{R})^{\frac{\gamma}{\nu}}\widehat{f}(t,\pi)\right\|^{2}_{\operatorname{HS}(\mathcal{H}_{\pi})}\mathrm{d}\mu(\pi)\nonumber\\
	&=&\|u^{f}(t,\cdot ; a,b)\|^{2}_{L^{2}_{\gamma+s\nu}(\mathbb{G})}+\|u^{f}(t,\cdot ; a,b)\|^{2}_{L^{2}_{\gamma}(\mathbb{G})}+\|f(t,\cdot)\|^{2}_{L^{2}_{\gamma}(\mathbb{G})}\nonumber\\
	&\lesssim&\|f\|^{2}_{C([0,T];L^{2}_{\gamma+s\nu}(\mathbb{G}))}+2\|f\|^{2}_{C([0,T];L^{2}_{\gamma}(\mathbb{G}))}\quad (\text{using } \eqref{eqnn2})\nonumber\\
	&\lesssim&\|f\|^{2}_{C([0,T];L^{2}_{\gamma+s\nu}(\mathbb{G}))}.
\end{eqnarray*}
This concludes the regularity estimate \eqref{NMI04}  and hence the proof of Lemma \ref{fnormm}.
\end{proof}
We have therefore determined the estimates for the non-homogeneous Cauchy problem  \eqref{homeq2}. 
We are now able to determine the regularity estimates for the Cauchy problem \eqref{mainpdenonhom}. 
    \begin{proof}[Proof of Theorem \ref{mthm2}](Regularity estimates) Recalling the relation \eqref{eqnnn} and utilizing the estimate \eqref{MI01} for  the solution $u^{u_{0}}(t,x)$ and  the estimate \eqref{NMI01} for the solution $u^{f}(t,x)$, we have
    	\begin{equation*}
    	\|u(t,\cdot)\|_{L^{2}(\mathbb{G})}\leq \left\|u^{u_{0}}(t,\cdot)\right\|_{L^{2}(\mathbb{G})}+\left\|u^{f}(t,\cdot)\right\|_{L^{2}(\mathbb{G})} 
    		\lesssim \|u_{0}\|_{L^{2}(\mathbb{G})}+\|f\|_{C([0,T];L^{2}(\mathbb{G}))}.
    	\end{equation*}
    	 This proves the regularity estimate \eqref{I01}. Similarily utilizing the estimate \eqref{MI02} for the solution $u^{u_{0}}(t,x)$ and the estimate \eqref{NMI02} for  the solution $u^{f}(t,x)$, we get
    	\begin{equation*}
    		\|u(t,\cdot)\|_{L^{2}_{\gamma}(\mathbb{G})}\leq \left\|u^{u_{0}}(t,\cdot)\right\|_{L^{2}_{\gamma}(\mathbb{G})}+\left\|u^{f}(t,\cdot)\right\|_{L^{2}_{\gamma}(\mathbb{G})} 
    		\lesssim \|u_{0}\|_{L^{2}_{\gamma}(\mathbb{G})}+\|f\|_{C([0,T];L^{2}_{\gamma}(\mathbb{G}))}.
    	\end{equation*}
    	This proves the regularity estimate \eqref{I02}. Further, utilizing the estimate \eqref{MI03} for the solution $u^{u_{0}}(t,x)$ and the estimate \eqref{NMI03} for  the solution $u^{f}(t,x)$, we can obtain
    	\begin{multline*}
    		\|\mathbb{D}_{(g)}u(t,\cdot)\|_{L^{2}(\mathbb{G})}\leq \left\|\mathbb{D}_{(g)}u^{u_{0}}(t,\cdot)\right\|_{L^{2}(\mathbb{G})}+\left\|\mathbb{D}_{(g)}u^{f}(t,\cdot)\right\|_{L^{2}(\mathbb{G})} 
    		\lesssim \|u_{0}\|_{L^{2}_{s\nu}(\mathbb{G})}+\\\|f\|_{C([0,T];L^{2}_{s\nu}(\mathbb{G}))}.
    	\end{multline*}
    	This proves the regularity estimate \eqref{I03}. Finally,  utilizing the estimate \eqref{MI04} for the solution $u^{u_{0}}(t,x)$ and the estimate \eqref{NMI04} for  the solution $u^{f}(t,x)$, we have
    	\begin{multline*}
    	\|\mathbb{D}_{(g)}u(t,\cdot)\|_{L^{2}_{\gamma}(\mathbb{G})}\leq \left\|\mathbb{D}_{(g)}u^{u_{0}}(t,\cdot)\right\|_{L^{2}_{\gamma}(\mathbb{G})}+\left\|\mathbb{D}_{(g)}u^{f}(t,\cdot)\right\|_{L^{2}_{\gamma}(\mathbb{G})} 
    	\lesssim\\ \|u_{0}\|_{L^{2}_{\gamma+s\nu}(\mathbb{G})}+\|f\|_{C([0,T];L^{2}_{\gamma+s\nu}(\mathbb{G}))}.
    \end{multline*}
    	This proves the  regularity estimate \eqref{I04} and hence complete the proof of Theorem \ref{mthm2}.   
\end{proof}
\section{Remarks}\label{remark}
We wrap up the article with a last comments about the assumption on coefficients in Theorem \ref{mthm1} and \ref{mthm2}:
\begin{enumerate}
	\item If we have the following additional information about the Rockland operator $\mathcal{R}$:
	\begin{equation*}
	\inf \pi(\mathcal{R}):=\inf\{\pi_{m}:m\in\mathbb{N}\}> 0, \text{ where } \pi(\mathcal{R})=\operatorname{diag}(\pi_{1}^{2},\pi_{2}^{2},\dots),
	\end{equation*}
	the assumption $0\leq a(t),b(t)\in C([0,T])$ such that $\inf\limits_{t\in[0,T]} b(t)=b_{0}>0$ on coefficients  in Theorem \ref{mthm2} can be relaxed to the assumption  $0\leq a(t),b(t)\in C([0,T])$ such that either  $\inf\limits_{t\in[0,T]} a(t)=a_{0}>0$ or  $\inf\limits_{t\in[0,T]} b(t)=b_{0}>0$.
	\item The above relaxation can be justified from the last inequality obtained in the estimate \eqref{festi}.
	\item In conclusion,  under these new assumptions on Rockland operator $\mathcal{R}$, the estimate of homogeneous Cauchy problem \eqref{mainpde} can be obtained from the estimates of non-homogeneous Cauchy problem \eqref{mainpdenonhom} by directly putting $f\equiv 0$.
\end{enumerate}
\section{Acknowledgement}
The first author was supported by Core Research Grant, RP03890G,  Science and Engineering
Research Board (SERB), DST,  India.   M. Ruzhansky is supported by the EPSRC Grants  EP/V005529/1 by the FWO Odysseus 1 grant G.0H94.18N: Analysis and Partial Differential Equations and by the  Methusalem programme of the Ghent University Special Research Fund (BOF) (Grantnumber
01M01021). The last author is supported by the institute assistantship from the Indian Institute of Technology Delhi, India.
\section{Declarations}
\noindent\textbf{Ethical Approval:} Not applicable.\\
 \textbf{Competing interests:} No potential competing of interest was reported by the author.\\
\textbf{Author's contributions:} 		 All the authors are contributed equally.\\
\textbf{Availability of data and materials:} 	The authors confirm that the data supporting the findings of this study are available within the article.
\bibliographystyle{alphaabbr}
\bibliography{time-fractional}

\end{document}